\documentclass[10pt,a4paper]{amsart}
\usepackage[utf8]{inputenc}
\usepackage{amsmath}
\usepackage{float}
\usepackage{geometry,graphicx,xcolor}
\usepackage{amssymb,color,esint,epic, graphicx,tikz}
\usepackage{parskip}
\usepackage{hyperref}
\usepackage{amsthm, amsfonts,enumerate}

\newtheorem{thm}{Theorem}[section]
\newtheorem{cor}[thm]{Corollary}

\newtheorem{prop}[thm]{Proposition}
\theoremstyle{definition}

\theoremstyle{remark}
\newtheorem{remark}{Remark}

\def\CC{\mathbb{C}}
\def\SS{\mathbb{S}^2}
\def\S^d{\mathbb{S}^d}
\def\E{\mathbb{E}}

\def\En{\mathcal{E}}

\newcommand{\per}{\text{per}}

\begin{document}

\title{Expected energy of zeros of elliptic polynomials}

\author{Victor de la Torre}
\address{Departament de Matem\`atiques i Inform\`atica, Universitat de Barcelona, Gran Via, 585, 08007 Barcelona, Spain}
\email{\href{mailto:delatorre@ub.edu }{\texttt{delatorre@ub.edu }}}
\author{Jordi Marzo}
\address{Departament de Matem\`atiques i Inform\`atica, Universitat de Barcelona, Gran Via, 585, 08007 Barcelona, Spain \vskip 0.01cm
Centre de Recerca Matemàtica, Edifici C, Campus Bellaterra, 08193 Bellaterra, Spain}
\email{\href{mailto:jmarzo@ub.edu}{\texttt{jmarzo@ub.edu}}}

\thanks{The authors would like to thank Joaquim Ortega-Cerd\`a and Carlos Beltr\'an for some helpful conversations.
\\
The authors have been partially supported by grant MTM2017-83499-P by the Ministerio de Ciencia, Innovaci\'on y Universidades, Gobierno de Espa\~na and by the 
Generalitat de Catalunya (project 2017 SGR 358).}

\begin{abstract}
    In 2011, Armentano, Beltr\'an and Shub obtained in \cite{ABS11} a closed expression for the expected logarithmic energy
     of the random point process on the sphere given by the roots of random elliptic polynomials. We consider a different approach which allows us to extend the study to the Riesz energies and to compute the expected separation distance. 
\end{abstract}

\date{\today}
\maketitle

\section{Introduction and Main results}

Elliptic polynomials, also called Kostlan-Shub-Smale or $SU(2)$ polynomials, are defined by 
\begin{equation*}
\sum_{n=0}^N a_n\sqrt{\binom{N}{n}}z^n,
\end{equation*}
where $a_n$ are i.i.d. random variables with standard complex Gaussian distribution. These polynomials appeared first in the mathematical physics literature \cite{BBL92,BBL96,Han96} and were readily 
studied from a mathematical point of view \cite{Kos93,SS93}. One reason for the interest in these polynomials is that the random point process on $\SS$ given by the stereographic 
projection of the roots of elliptic polynomials is invariant through rotations. Moreover, it is the unique point process given by zeros of random analytic functions with this property \cite{Sod00}. Among its many interesting properties, especially 
relevant are the connections, studied in \cite{SS93}, with well conditioned polynomials and with minimal logarithmic energy points.

The Riesz or logarithmic energy of a set of $N$ different points $x_1,\ldots,x_N$
on the unit sphere $\SS \subset \mathbb R^3$ is 
$$E_s(x_1,\dots ,x_N)=\sum_{i\neq j}f_s(|x_i-x_j|),$$ 
where $f_s(r)=r^{-s}$ for $s\neq 0$ and $f_0(r)=-\log r$ are, respectively, the Riesz and logarithmic potentials. We denote 
the extremal (minimal or maximal) energy attained by a set of $N$ points on the sphere by
$$\mathcal{E}_s(N)=
\left\{
	\begin{array}{ll}
		\min_{x_1,\dots , x_N\in \SS} E_s(x_1,\dots ,x_N)  & \mbox{if } s\ge 0, \\
		\max_{x_1,\dots , x_N\in \SS} E_s(x_1,\dots ,x_N) & \mbox{if } s < 0.
	\end{array}
\right.
$$ 

The condition number of a univariate polynomial, defined by Shub and Smale, is a measure of how much the roots of a polynomial change when perturbing the coefficients. 
It was shown in \cite{SS93b} that points of almost minimal logarithmic energy, $s=0,$ are the roots of well conditioned polynomials. In \cite{SS93}, the authors also proved 
that, with high probability, elliptic polynomials are well conditioned, see  \cite{BEMOC21,BL20} for a deterministic example. 
It was therefore natural to study the expected energy of the zeros of elliptic polynomials. This was done in \cite{ABS11}, where the authors
obtained the following closed expression for the expected logarithmic energy of random points $x_1,\dots , x_N\in \mathbb S^2,$ 
images by the stereographic projection of zeros of elliptic polynomials,
\begin{equation} \label{armentano}
\mathbb E [ E_0(x_1,\dots , x_N)]=\left(\frac{1}{2}-\log 2\right)N^2
-\frac{1}{2}N\log N
-\left(\frac{1}{2}-\log 2\right)N.
\end{equation}
The asymptotic expression above is indeed very close to the minimal logarithmic energy of $N$ points on the sphere, see Section \ref{comparison_results}. Working in a more general setting, in \cite{Zho08, ZZ10}
the same expression (\ref{armentano}) was obtained but with a $o(N)$ remainder. Our main result is an extension of the above result (\ref{armentano}) to the Riesz $s$-energies for $s<4.$


\begin{thm}\label{asymptotics}
Let $x_1,\dots, x_N\in \mathbb S^2$ be the image by the stereographic projection of $N$ points drawn from zeros of elliptic polynomials. Then, 
\begin{enumerate}[(i)]
\item
for $s<4,$ $s\neq 0,2$ and a fixed $m\geq 1$,
\begin{multline}\label{newrieszasymptotics}
\mathbb E [ E_s(x_1,\dots , x_N)]=\frac{2^{1-s}}{2-s}N^2\\
+\frac{\Gamma\left(1-\frac{s}{2}\right)}{2^{s+1}}\Bigg[s\left(1+\frac{s}{2}\right)\sum_{j=0}^{m-1}\frac{B_{2j}^{(\frac{s}{2})}(\frac{s}{4})(1-\frac{s}{2})_{2j}}{(2j)!}N^{\frac{s}{2}+1-2j}\zeta\left(1-\frac{s}{2}+2j,1+\frac{4-s}{4N}\right)
\\
+s\left(1-\frac{s}{2}\right)\sum_{j=0}^{m-1}\frac{B_{2j}^{(\frac{s}{2}-1)}(\frac{s-2}{4})(2-\frac{s}{2})_{2j}}{(2j)!}N^{\frac{s}{2}-2j}\zeta\left(2-\frac{s}{2}+2j,1+\frac{2-s}{4N}\right)\Bigg]
+ O\left(N^{\frac{s}{2}+1-2m}\right),
\end{multline}
for $N\to +\infty.$

\item Moreover, the energies with $s=-2n$ for integer $n\geq -1$ can be computed exactly:

For $s=0$,
\begin{equation}\label{ourABS}
 \mathbb E [ E_0(x_1,\dots , x_N)]=\left(\frac{1}{2}-\log 2\right)N^2
-\frac{N\log N}{2}
-\left(\frac{1}{2}-\log 2\right)N.   
\end{equation}

For $s=2$,
\begin{equation}\label{2energy}
    \mathbb E [ E_2(x_1,\dots , x_N)]=-\frac{ N\pi }{4}\sum_{j=1}^{N-1}\frac{j}{N}\cot\left(\frac{\pi j}{N}\right)+\frac{3N^2}{8}-\frac{3N}{8}.
\end{equation}

For $s=-2n$, $n\geq 1$,
\begin{align}\label{evenenergies}
 \mathbb E  & [ E_{-2n}(x_1,\dots , x_N)]
=2^{2 n} N^2 \left(\frac{1}{n + 1} - \frac{n (n - 1)}{n + 1} - 
    n \sum_{m=1}^{n+1}\frac{1}{m}\right)
    \\
    &
		+  2^{2 n}n N\left(-\gamma + 
    \sum_{m=1}^{n+1}\binom{n + 1}{m}(-1)^m \psi\left(\frac{m}{N}\right)\left(\frac{n - 1 }{n + 1}m + 1\right)\right). \nonumber
\end{align}
\end{enumerate}
\end{thm}


In the above result, $\gamma$ is the Euler-Mascheroni constant, $B_{2j}^{(2\rho)}(\rho)$ are the generalized Bernoulli polynomials defined by
$$\left(\frac{t}{e^{t}-1}\right)^{2\rho}e^{\rho t}=\sum_{j=0}^\infty \frac{t^{2j}}{(2j)!}B_{2j}^{(2\rho)}(\rho),$$
for $|t|<2\pi$, with $B_{0}^{(2\rho)}(\rho)=1,$  
$$\zeta (s,a)=\sum_{j=0}^\infty \frac{1}{(j+a)^s},\;\;\Re s>1, \; a\not\in \mathbb Z_{\le 0}$$ is the Hurwitz Zeta function and $\psi(z)=\Gamma'(z)/\Gamma(z)$ is the digamma function.


By considering two terms of the asymptotic expansion of the Hurwitz Zeta function
\[
\zeta(s,1+a)=\sum_{k=0}^{\infty} \frac{(-1)^k(s)_k\zeta(s+k)}{k!}a^k,
\]
for $|a|<1$ and $s\neq 1$ \cite[25.11.10]{DLMF} and taking $m=1$ in \eqref{newrieszasymptotics} we get, for $0,2\neq s<4$,
\begin{equation}\label{rieszasymptotics}
\E[E_s(x_1,\dots , x_N)]=\frac{2^{1-s}}{2-s}N^2+C(s)N^{1+s/2}+\frac{s}{16}C(s-2) N^{s/2}+O(N^{-1+s/2}),
\end{equation}
when $N\to \infty,$ where
\begin{gather}
C(s)=\frac{1}{2^s}\frac{s}{2}\left(1+\frac{s}{2}\right)\Gamma\left(1-\frac{s}{2}\right)\zeta\left(1-\frac{s}{2}\right) \label{constant1}.
\end{gather}

\begin{remark}
The result above for the expected Riesz energy allows us to compare the zeros of elliptic polynomials with other point processes, for example in terms of expected $p$-moments of averages. Indeed, from Khintchine's 
inequality \cite[Theorem 3]{KK01}, it follows that 
$$\E\left[\left| \sum_{i=1}^N x_i\right|^p\right]\sim N^{p/2}$$ 
when $x_1,\dots, x_N$ are uniform i.i.d. points on the sphere $\mathbb S^2$ and $1\le p<\infty.$ For points drawn from the spherical ensemble,
for which there is repulsion between points, it follows from 
\begin{equation}\label{equation-2}
\sum_{i,j=1}^N |x_i-x_j|^{2}=2N^2-2\left|\sum_{i=1}^N x_i\right|^2
\end{equation}
and the results about the expected Riesz energy $s=-2$ in \cite{AZ15} that the expected 2-moment is bounded. Hence, for the spherical ensemble $\E\left[\left| \sum_{i=1}^N x_i\right|^p\right]$ is bounded 
for $1\le p\le 2,$ and numerical 
simulations suggest that the same holds for $p>2.$ 
In our case, for zeros of elliptic polynomials mapped to the sphere by the stereographic projection, it follows from (\ref{evenenergies}) that
\begin{equation}
\E\left[\left| \sum_{i=1}^N x_i\right|^2\right]=4\frac{\zeta(3)}{N}+o (N^{-1} ),
\end{equation}
for $N\to +\infty,$
and the average $p$-moments for $1\le p\le 2$ converge to zero. Again, numerical simulations suggest the same behavior for $p>2$. 
It is well known that minimal logarithmic points have center of mass in the center of the sphere, i.e. have zero dipole, \cite[Corollary 6.7.5]{BHS19}, \cite{BBP94}. Therefore, the behavior of the expected $p$-moments matches the particularly 
low logarithmic energy of zeros of elliptic polynomials.
For the comparison with minimal and expected energies of other point processes, 
see discussion in Section \ref{comparison_results}.
\end{remark}

\begin{figure}[H]
\centering
       \includegraphics[width=0.6\textwidth]{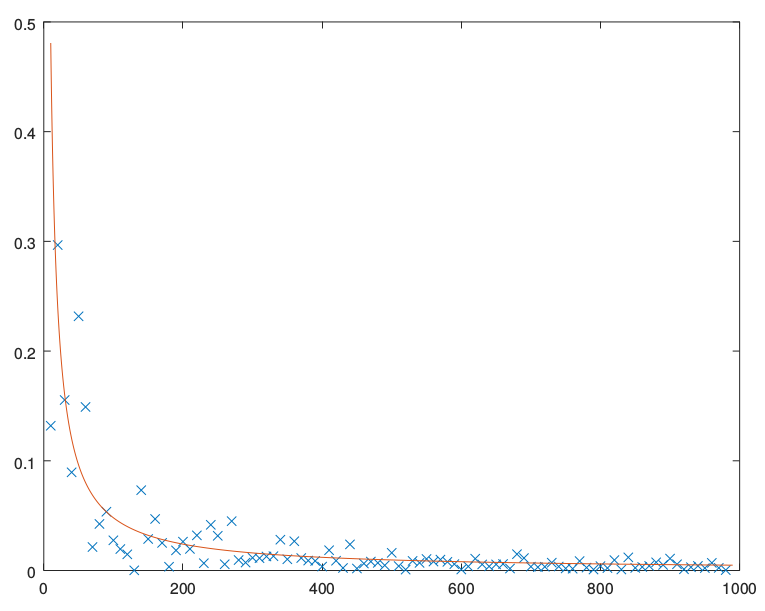}
    \caption{Plot of $4\zeta(3)/N$ and realizations of $| \sum_{i=1}^N x_i |^2$ for natural $N$ up to 1000.}
    \label{fig:center_of_mass}
\end{figure}

%
%
%


In our last result, we compute a closed expression for the expected separation distance between 
points drawn from zeros of elliptic polynomials. 
The separation distance of $X_N=\{ x_1,\dots , x_N \}$ is defined by
$$\mbox{sep}(X_N)=\min_{i\neq j}|x_i-x_j |,$$
and its counting version by
$G(t,X_N)=\sharp\{i<j:|x_i-x_j|\leq t\}.$ Recall that energy minimizers have a separation distance of order $N^{-1/2},$ \cite[Section 6.9]{BHS19}.

\begin{thm} \label{separation_thm}
Let $X_N$ be a set of $N-$points drawn from zeros of elliptic polynomials mapped to the sphere by the stereographic projection. Then
\begin{equation}\label{the_sepa_result}
 \E[G(t,X_N)]=\frac{t^2 N^2 }{8}
-\frac{N}{2}
+\frac{t^2 N^2}{8(4-t^2)\left(\left(\frac{4}{4-t^2}\right)^N-1\right)}\left[8-t^2-t^2N-\frac{t^2N}{\left(\frac{4}{4-t^2}\right)^N-1}\right] .
\end{equation}
Therefore,
\begin{equation}
    \label{sep_equality}
    \E[G(t,X_N)]=\frac{N^3 t^4}{128}(1+o(1)),
\end{equation}
if $t=o(1/\sqrt{N})$, and moreover
\begin{equation}
    \label{sep_inequality}
    \E[G(t,X_N)]\leq \frac{N^3 t^4}{128},
\end{equation}
for $t\le 2.$
\end{thm}

\begin{figure}[H]
\centering
       \includegraphics[width=0.6\textwidth]{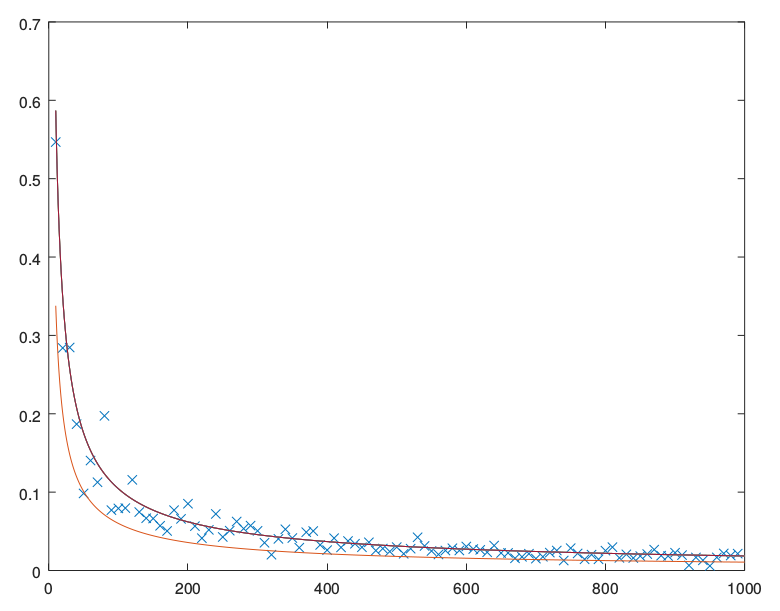}
    \caption{x marks correspond to the values of the minimal separation for realizations of $N$ elliptic zeros for natural $N$ from 10 up to 1000. The continuous graph are $c N^{-3/4}$ for $c=1.89$ (yellow) and $3.27$ (brown): using Chebyshev's inequality at least $90\%$ of the
    realizations are above yellow and at least $10\%$ above brown.}
    \label{fig:separation}
\end{figure}

Note that $\mbox{sep}(X_N)\le t$
implies $G ( t ,X_N )\ge 1$, hence 
$\mathbb P ( \mbox{sep} ( X_N )\le t )\le \mathbb  P ( G ( t , X_N )\ge 1 )\le  \mathbb E ( G ( t , X_N ))$ and therefore, as in the harmonic case, see \cite{BMOC16}, an $N$-tuple drawn from the zeros of elliptic polynomials likely
satisfies $\mbox{sep}(X_N)=\Omega(N^{-3/4}),$ Figure \ref{fig:separation}. See also \cite[Corollary 1.6]{AZ15} for the analogue result for the spherical ensemble. 

\subsection{Organization of the paper}
In section \ref{sec_intensity} we compute the 2-point intensity function of our point process and explain how to compute the expected energy.
Section \ref{sec_proof_main} contains the proof of our main result, Theorem \ref{asymptotics}. In section 
\ref{comparison_results} we deduce some bounds for the extremal energy and compare our bounds with previous results.
Finally, in section \ref{sec_separation} we prove Theorem \ref{separation_thm} about separation.



\section{Intensity function and expected Riesz energy}  \label{sec_intensity}

In this section we compute the $2$-point intensity function of the random point process on $\SS$ corresponding to the stereographic projection of the roots of random elliptic polynomials
\begin{equation*}
P_N(z)=\sum_{n=0}^N a_n\sqrt{\binom{N}{n}}z^n,
\end{equation*}
where $a_n$ are i.i.d. random variables with standard complex Gaussian distribution.
Let $F(x,y)$ be a
measurable function defined on $\SS\times \SS$ whose variables will be considered in $\CC$ through the stereographic projection, i.e. $F(z,w)=F(x(z),y(w))$, with the points $x,y\in\SS$ corresponding to $z,w\in\CC.$ By Campbell's formula, 
if $x_1,\dots x_N\in \mathbb S^2$ are the images of the zeros $z_1,\dots, z_N$ of elliptic polynomials, then
\begin{equation}\label{computation1}
\begin{aligned}
\E & \left[\sum_{i\neq j}F(x_{i},x_{j})\right]&=\E\left[\sum_{i\neq j}F(z_i,z_j)\right]=\int_{\CC}\int_{\CC}F(z,w)\rho_2(z,w)dz dw,
\end{aligned}
\end{equation}
with $\rho_2(z,w)$ the 2-point intensity function given by \cite[Corollary 3.4.2]{HKPV09} 
\begin{equation}\label{rhoexpression}
\rho_{2}(z_{1},z_2)=\frac{\per(C-BA^{-1}B^{*})}{\det(\pi A)},
\end{equation}
where $A,B,C$ are the $2\times 2$ matrices
\begin{align*}
A(i,j) &= \E[P_N(z_{i})\overline{P_N}(z_{j})],\\
B(i,j) &= \E[P_N'(z_{i})\overline{P_N}(z_{j})],\\
C(i,j) &= \E[P_N'(z_{i})\overline{P_N'}(z_{j})].
\end{align*}

\begin{figure}
\centering
       \includegraphics[width=0.4\textwidth]{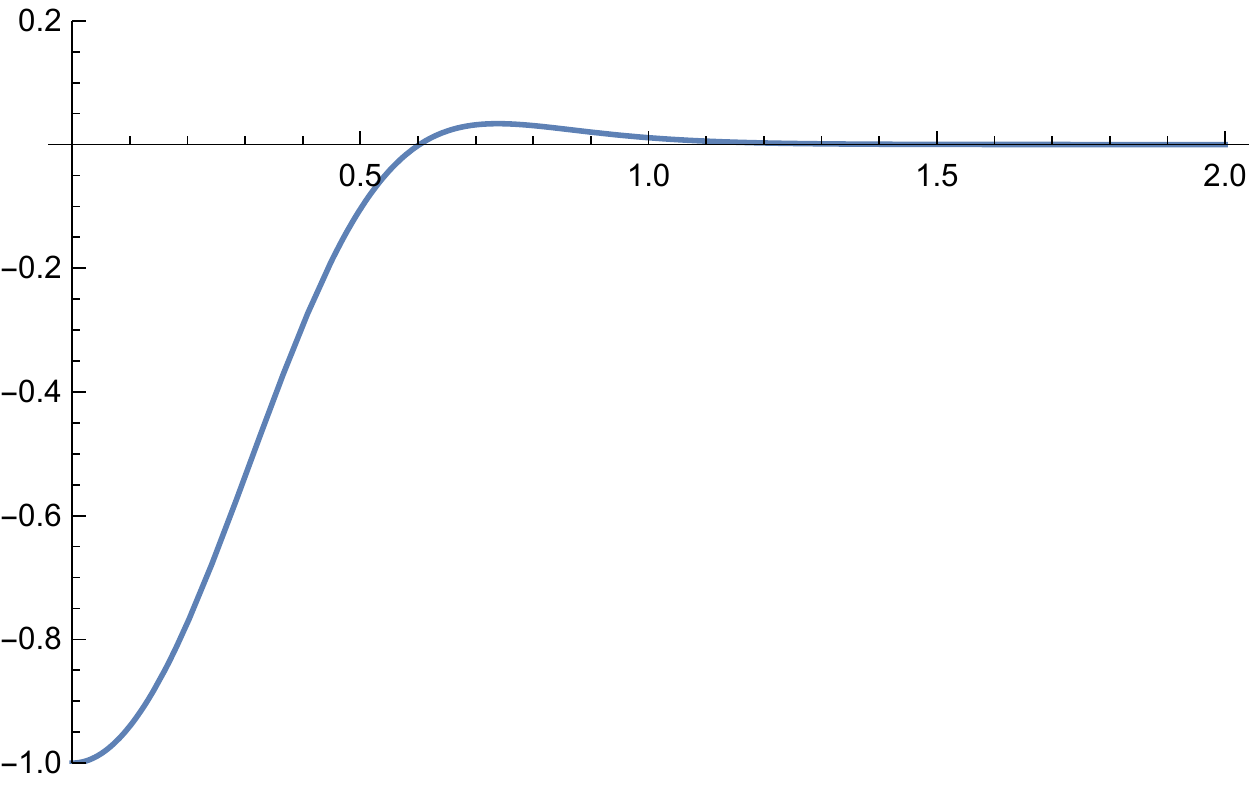}
    \caption{$\frac{\pi^2}{N^2}(\rho_2(r,0)-\rho_1(r)\rho_1(0))$ for $r>0$ and $N=10$.}
    \label{fig:det}
\end{figure}

It is easy to see that when $F$ is rotational invariant we get 
\begin{equation}\label{computation2}
\E\left[\sum_{i\neq j}F(x_{i},x_{j})\right]=\pi\int_{\CC}F(z,0)\rho_2(z,0)dz.
\end{equation}
Therefore, it is enough to compute $\rho_2(z_1,z_2)$ for $z_1=z\in \mathbb C$ and $z_2=0$. The matrices in \eqref{rhoexpression} are then
$$
A=
\begin{pmatrix}
(1+|z|^2)^N&1\\
1&1
\end{pmatrix},$$
$$B=
N\begin{pmatrix}
\overline{z}(1+|z|^2)^{N-1}&0\\
\overline{z}&0
\end{pmatrix},$$
$$C=
N\begin{pmatrix}
(1+|z|^2)^{N-2}(1+N|z|^2)&1\\
1&1
\end{pmatrix},
$$
and we obtain 
$$
 \rho_{2}(z,0)=\frac{N^2 \left[ \left(1-\frac{N|z|^2}{(1+|z|^2)^N-1}\right)^2(1+|z|^2)^{N-2}+\left(1-\frac{N|z|^2(1+|z|^2)^{N-1}}{(1+|z|^2)^N-1}\right)^2\right]}{\pi^2[(1+|z|^2)^N-1]},
$$
see \cite{Han96} and Figure \ref{fig:det} where one can notice that this point process is not determinantal (\cite[p.83]{HKPV09}).

Using the relation with the chordal metric
\[|x-y|=\frac{2|z-w|}{\sqrt{1+|z|^2}\sqrt{1+|w|^2}},\]
we get
\begin{equation*}
\begin{aligned}
\E[E_0(x_1,\dots x_N)]&=\pi\int_{\CC}f_0 (z,0)\rho_2(z,0)dz=-\pi\int_{\CC}\log\left(\frac{2|z|}{\sqrt{1+|z|^2}}\right)\rho_2(z,0)dz\\
&=-2N^2\int_0^\infty r\log\left(\frac{2r}{\sqrt{1+r^2}}\right)\gamma(r)dr,
\end{aligned}
\end{equation*}
and for $s\neq 0$
\begin{equation}\label{sRiesz}
\begin{aligned}
\E[E_{s}(x_1,\dots , x_N)]&=\pi\int_{\CC}f_s(z,0)\rho_2(z,0)dz=\pi\int_{\CC}\left(\frac{2|z|}{\sqrt{1+|z|^2}}\right)^{-s}\rho_2(z,0)dz\\
&=2^{1-s}N^2\int_0^\infty r^{1-s}(1+r^2)^{s/2}\gamma(r)dr,
\end{aligned}
\end{equation}
where
\begin{equation}\label{funciogamma}
\gamma(r)=\frac{\left[ \left(1-\frac{Nr^2}{(1+r^2)^N-1}\right)^2(1+r^2)^{N-2}+\left(1-\frac{Nr^2(1+r^2)^{N-1}}{(1+r^2)^N-1}\right)^2\right]}{[(1+r^2)^N-1]}.
\end{equation}
In the logarithmic case, one can compute directly a primitive function that leads to the correct energy \eqref{armentano}. However, we will compute the expected logarithmic energy as the limit of the Riesz case at $s=0$.


\section{Proof of theorem \ref{asymptotics}}\label{sec_proof_main}

In this section we prove first our general result (\ref{newrieszasymptotics}) with the auxiliary Proposition \ref{aux_prop}. Then we prove the cases (\ref{2energy}),(\ref{evenenergies}) and finally (\ref{ourABS}).

\begin{proof}
To simplify the notation we write $\E[E_{s}]$ instead of $\E[E_{s}(x_1,\dots , x_N)].$
The change of variables $r=\sqrt{x}$ in \eqref{sRiesz} yields
\begin{multline*}
\E[E_{s}]=\frac{N^2}{2^s}\int_0^\infty \frac{x^{-s/2}(1+x)^{s/2}}{[(1+x)^N-1]^3}\\
\left[\left((1+x)^N-1-Nx\right)^2(1+x)^{N-2}+\left((1+x)^N-1-Nx(1+x)^{N-1}\right)^2\right]dx.
\end{multline*}
The integrand is equivalent to $x^{-2}$ at infinity, which is integrable, and to $x^{1-s/2}$ at $x=0$, which is integrable iff $1-s/2>-1$. Then, the energy will be finite iff $s<4$.

Now let us compute the integral. We take $r=s/2$ for simplicity, so we will be assuming $r<2$ throughout the proof. Using that $\frac{1}{(x-1)^3}=\frac{1}{2}\sum_{k=2}^\infty k(k-1)x^{-(k+1)}$ for $x>1$ and the fact that all the terms are positive, we get
\begin{align}\label{begin_s}
& \E[E_{2r}]=\frac{N^2}{2^{2r+1}}\sum_{k=2}^\infty k(k-1)\int_0^\infty \frac{(1+x)^{r-N(k+1)}}{x^r}
\\
&
\left[\left((1+x)^N-1-Nx\right)^2(1+x)^{N-2}+\left((1+x)^N-1-Nx(1+x)^{N-1}\right)^2\right]dx
\nonumber
\\
&
=
\frac{N^2}{2^{2r+1}}\lim_{M\to\infty}\sum_{k=2}^M k(k-1)\left[\underbrace{\int_0^\infty \frac{\left[(1+x)^{r-2-Nk}+(1+x)^{r-N(k+1)}\right]\left((1+x)^N-1\right)^2}{x^{r}} dx}_{A_k} \right.
\nonumber
\\
&
-2N\underbrace{\int_0^\infty x^{1-r}\left[(1+x)^{r-2-Nk}+(1+x)^{r-1-Nk}\right]\left((1+x)^N-1\right) dx}_{B_k}
\nonumber
\\
&
\left.
+N^2\underbrace{\int_0^\infty x^{2-r}\left[(1+x)^{r-2-Nk}+(1+x)^{r-2-N(k-1)}\right] dx}_{C_k}\right].
\nonumber
\end{align}
Using the following integral representation for the beta function (see \cite[8.380 (3)]{GR07}),
\[
B(x,y)=\frac{\Gamma(x)\Gamma(y)}{\Gamma(x+y)}=\int_0^\infty \frac{t^{x-1}}{(1+t)^{x+y}}\, dt\qquad x,y>0,
\]
it is immediate to obtain $B_k,C_k$ in (\ref{begin_s})
\begin{equation}\label{Bk}
B_k=B(2-r,N(k-1))-B(2-r,Nk)+B(2-r,N(k-1)-1)-B(2-r,Nk-1),    
\end{equation}
\begin{equation}\label{Ck}
C_k=B(3-r,Nk-1)+B(3-r,N(k-1)-1),    
\end{equation}
so
\begin{align*}
& -2NB_k+N^2C_k 
\\
&
=\Gamma(2-r)\left[-2N\left(\frac{\Gamma(N(k-1))}{\Gamma(N(k-1)+2-r)}-\frac{\Gamma(Nk)}{\Gamma(Nk+2-r)}+\frac{\Gamma(N(k-1)-1)}{\Gamma(N(k-1)+1-r)} \right.\right.
\\
&
\left.\left.
-\frac{\Gamma(Nk-1)}{\Gamma(Nk+1-r)}\right)
+N^2(2-r)\left(\frac{\Gamma(Nk-1)}{\Gamma(Nk+2-r)}+\frac{\Gamma(N(k-1)-1)}{\Gamma(N(k-1)+2-r)}\right)\right].
\end{align*}
To compute $A_k$ we integrate by parts. Let $\beta\in\{r-2-Nk,r-N(k+1)\}$ denote the exponent in $(1+x)$. If $r\neq 1$,
\begin{align*}
\int_0^\infty & \frac{(1+x)^\beta \left((1+x)^N-1\right)^2}{x^{r}} dx=\frac{1}{1-r}x^{1-r}(1+x)^\beta\left((1+x)^N-1\right)^2\Big|_0^\infty
\\
- & \frac{1}{1-r}\int_0^\infty x^{1-r}\left[\beta(1+x)^{\beta-1}\left((1+x)^N-1\right)^2+2N(1+x)^{\beta+N-1} \left((1+x)^N-1\right)\right]dx
\\
=
&
\frac{-1}{1-r}\left[\beta\int_0^\infty x^{1-r}(1+x)^{\beta-1}\left((1+x)^{2N}-2(1+x)^N+1\right)dx \right.
\\
+
&
\left.
2N\int_0^\infty x^{1-r}(1+x)^{\beta+N-1}\left((1+x)^N-1\right)  dx\right]
\\
=
&
\frac{-1}{1-r}\left[\beta\left(B(2-r,-\beta-2N-1+r)-2B(2-r,-\beta-N-1+r)+B(2-r,-\beta-1+r)\right) \right. 
\\
+
&
\left.
2N\left(B(2-r,-\beta-2N-1+r)-B(2-r,-\beta-N-1+r)\right)\right]
\\
=
&
B(1-r,-\beta-2N-1+r)-2B(1-r,-\beta-N-1+r)+B(1-r,-\beta-1+r).
\end{align*}
%
Then
\begin{align}\label{Ak}
A_k & =
B(1-r,N(k-2)+1)-2B(1-r,N(k-1)+1)+B(1-r,Nk+1)
\\
&
+B(1-r,N(k-1)-1)-2B(1-r,Nk-1)+B(1-r,N(k+1)-1),\nonumber
\end{align}
or, in terms of the gamma function,
\begin{align*}
A_k & =\Gamma(1-r)\left[\frac{\Gamma(N(k-2)+1)}{\Gamma(N(k-2)+2-r)}-2\frac{\Gamma(N(k-1)+1)}{\Gamma(N(k-1)+2-r)}+\frac{\Gamma(Nk+1)}{\Gamma(Nk+2-r)} \right.
\\
&
\left.
+\frac{\Gamma(N(k-1)-1)}{\Gamma(N(k-1)-r)}-2\frac{\Gamma(Nk-1)}{\Gamma(Nk-r)}+\frac{\Gamma(N(k+1)-1)}{\Gamma(N(k+1)-r)}\right],
\end{align*}
provided that $r\neq 1$. The case $r=1$ will be studied as the limit $r\to 1$.

Therefore, for $r\neq 1$, writing all together
\begin{align}\label{take_limit}
&  \E[E_{2r}] =\frac{N^2}{2^{2r+1}}\lim_{M\to\infty}\left[\sum_{k=2}^M k(k-1)\Gamma(1-r)\left(\frac{\Gamma(N(k-2)+1)}{\Gamma(N(k-2)+2-r)}-2\frac{\Gamma(N(k-1)+1)}{\Gamma(N(k-1)+2-r)}\right.\right.
\\
&
\left.+\frac{\Gamma(Nk+1)}{\Gamma(Nk+2-r)}+\frac{\Gamma(N(k-1)-1)}{\Gamma(N(k-1)-r)}-2\frac{\Gamma(Nk-1)}{\Gamma(Nk-r)}+\frac{\Gamma(N(k+1)-1)}{\Gamma(N(k+1)-r)}\right) \nonumber
\\
&
+\sum_{k=2}^M k(k-1)\Gamma(2-r)\left(-2N\left(\frac{\Gamma(N(k-1))}{\Gamma(N(k-1)+2-r)}-\frac{\Gamma(Nk)}{\Gamma(Nk+2-r)} +\frac{\Gamma(N(k-1)-1)}{\Gamma(N(k-1)+1-r)} \right.\right.  \nonumber
\\
&
\left.\left.\left.
-\frac{\Gamma(Nk-1)}{\Gamma(Nk+1-r)}\right)+N^2(2-r)\left(\frac{\Gamma(Nk-1)}{\Gamma(Nk+2-r)}+\frac{\Gamma(N(k-1)-1)}{\Gamma(N(k-1)+2-r)}\right)\right)\right]. \nonumber
\end{align}

The sums get simplified by using the property $\Gamma(z+1)=z\Gamma(z)$ and changing the indices in such a way that all quotients have the form $\Gamma(Nk+1)/\Gamma(Nk+2-r)$
\begin{align*}
& \E[E_{2r}]=\frac{\Gamma(1-r) N^2}{2^{2r+1}}
\\
&
\lim_{M\to\infty}\left[\frac{2}{\Gamma(2-r)}+\sum_{k=1}^M (1-r+Nk(1+r))\frac{2r \Gamma(Nk)}{\Gamma(Nk+2-r)}-(M+1)M\frac{\Gamma(N(M-1)+1)}{\Gamma(N(M-1)+2-r)} \right.
\\
&
-\frac{(M+1) \left(r \left(N \left(N (4 M-3)-2 M+2\right)-1\right)-2 (N-1) N (M-1)+(N-1)^2 r^2\right)}{N(NM-1)}
\\
&
\left.
\frac{\Gamma(NM+1)}{\Gamma(NM+2-r)}+\frac{M(M-1)(N(M+1)-r)(N(M+1)+1-r)}{N(M+1)(N(M+1)-1)}\frac{\Gamma(N(M+1)+1)}{\Gamma(N(M+1)+2-r)}\right].
\end{align*}
%
Taking the asymptotic expansion of the terms in $M$ as $M\to \infty$, we get
\begin{align}\label{expression_M}
& \E[E_{2r}]=\frac{\Gamma(1-r)N^2}{2^{2r+1}}\lim_{M\to\infty}\left[\frac{2}{\Gamma(2-r)}+2r(1+r)\sum_{k=1}^M \frac{\Gamma(Nk+1)}{\Gamma(Nk+2-r)} \right.
\\
&
\left.
+2r(1-r)\sum_{k=1}^M \frac{\Gamma(Nk)}{\Gamma(Nk+2-r)}
-2(1+r)N^{r-1}M^r-r(N+r+Nr-r^2)N^{r-2}M^{r-1} \right]. \nonumber
\end{align}
Applying Proposition \ref{aux_prop} below we obtain the following expression for every $r\neq 0,1$ with $r<2$
\begin{align*}
& \E[E_{2r}]=\frac{\Gamma(1-r)N^2}{2^{2r+1}}
\\
&
\left[\frac{2}{\Gamma(2-r)}+2r(1+r)\sum_{j=0}^{m-1}\frac{B_{2j}^{(r)}(\frac{r}{2})(1-r)_{2j}}{(2j)!}\zeta\left(1-r+2j,1+\frac{2-r}{2N}\right) N^{r-1-2j} \right.
\\
&
\left.
+2r(1-r)\sum_{j=0}^{m-1}\frac{B_{2j}^{(r-1)}(\frac{r-1}{2})(2-r)_{2j}}{(2j)!}\zeta\left(2-r+2j,1+\frac{1-r}{2N}\right)N^{r-2-2j}+ O\left(N^{r-1-2m}\right)\right].
\end{align*}
Writing the expression in terms of $s=2r$ yields the result (\ref{newrieszasymptotics}).

Now we prove (\ref{2energy}), i.e. $r=1$, from the case $r\neq 1.$ By continuity, the evaluation of the integral at the beginning of \eqref{begin_s} can be performed by taking the limit $r\to 1$ in $A_k$, $B_k$, $C_k$, that is, in 
both sums in \eqref{take_limit}. The only tricky limit is the first one. It can be computed using the asymptotic expansion
\[
\frac{1}{\Gamma(a+\gamma)}=\frac{1}{\Gamma(a)}-\frac{\psi(a)}{\Gamma(a)}\gamma+o \left(\gamma\right),
\]
for $\gamma\to 0,$ where $a$ will be a natural number. Considering $\gamma = 1-r$,
\begin{align*}
\lim_{\gamma\to 0} & \Gamma(\gamma)  \Bigg[\frac{\Gamma(N(k-2)+1)}{\Gamma(N(k-2)+1+\gamma)}-\frac{2\Gamma(N(k-1)+1)}{\Gamma(N(k-1)+1+\gamma)}+\frac{\Gamma(Nk+1)}{\Gamma(Nk+1+\gamma)}+\frac{\Gamma(N(k-1)-1)}{\Gamma(N(k-1)-1+\gamma)}
\\
&
-2\frac{\Gamma(Nk-1)}{\Gamma(Nk-1+\gamma)}+\frac{\Gamma(N(k+1)-1)}{\Gamma(N(k+1)-1+\gamma)}\Bigg]=-\psi(N(k-2)+1)+2\psi(N(k-1)+1)
\\
&
-\psi(Nk+1)-\psi(N(k-1)-1)+2\psi(Nk-1)-\psi(N(k+1)-1),
\end{align*}
and we get from (\ref{take_limit})
\begin{align}\label{s2_twosums}
& \E[E_{2}]=\frac{N^2}{2^3}\lim_{M\to\infty}\left[\sum_{k=2}^M k(k-1)(-\psi(N(k-2)+1)+2\psi(N(k-1)+1) \right.
\\
&
-\psi(Nk+1) -\psi(N(k-1)-1)+2\psi(Nk-1)-\psi(N(k+1)-1)) \nonumber
\\
&
+\sum_{k=2}^M k(k-1)\left(-2N\left(\frac{\Gamma(N(k-1))}{\Gamma(N(k-1)+1)}-\frac{\Gamma(Nk)}{\Gamma(Nk+1)}+\frac{\Gamma(N(k-1)-1)}{\Gamma(N(k-1))} \right.\right.  \nonumber
\\
&
\left.\left.\left.
-\frac{\Gamma(Nk-1)}{\Gamma(Nk)}\right)+N^2\left(\frac{\Gamma(Nk-1)}{\Gamma(Nk+1)}+\frac{\Gamma(N(k-1)-1)}{\Gamma(N(k-1)+1)}\right)\right)\right]. \nonumber
\end{align}
The first sum in \eqref{s2_twosums} can be rewritten as
\begin{align*}
& \Sigma_1  :=\sum_{k=2}^M k(k-1)(-\psi(N(k-2)+1)+2\psi(N(k-1)+1)-\psi(Nk+1)
\\
&
-\psi(N(k-1)-1)+2\psi(Nk-1)-\psi(N(k+1)-1))
\\
&
=(M+2)(M+1)\psi(NM+1)+(M+1)M\psi(N(M-1)+1)
\\
&
-2(M+1)M \psi(NM+1)-2\sum_{k=0}^M \psi(Nk+1)
\\
&
+(M+1)M \psi(NM-1)-M(M-1)\psi(N(M+1)-1)-2\sum_{k=1}^M\psi(Nk-1),
\end{align*}
while the second becomes
\begin{align*}
& \Sigma_2  :=\sum_{k=2}^M k(k-1)\left(-2N\left(\frac{\Gamma(N(k-1))}{\Gamma(N(k-1)+1)}-\frac{\Gamma(Nk)}{\Gamma(Nk+1)}+\frac{\Gamma(N(k-1)-1)}{\Gamma(N(k-1))} \right.\right.
\\
&
\left.\left.
-\frac{\Gamma(Nk-1)}{\Gamma(Nk)}\right)+N^2\left(\frac{\Gamma(Nk-1)}{\Gamma(Nk+1)}+\frac{\Gamma(N(k-1)-1)}{\Gamma(N(k-1)+1)}\right)\right)
\\
&
=\sum_{k=2}^M k(k-1)\left(-2N\left(\frac{1}{N(k-1)}-\frac{1}{Nk}+\frac{1}{N(k-1)-1}-\frac{1}{Nk-1}\right)\right.
\\
&
+\left.N^2\left(\frac{1}{Nk-1}-\frac{1}{Nk}+\frac{1}{N(k-1)-1}-\frac{1}{N(k-1)}\right)\right)
\\
&
=\sum_{k=2}^{M}k(k-1)\left(\frac{2N-N^2}{Nk}-\frac{2N+N^2}{N(k-1)}+\frac{(2N+N^2)}{Nk-1}-\frac{(2N-N^2)}{N(k-1)-1}\right)
\\
&
=-\sum_{k=1}^{M-1}\frac{2kN(2+kN)}{Nk}+\frac{M(M-1)}{NM}(2N-N^2)+\sum_{k=1}^{M-1}\frac{2kN(-2+kN)}{Nk-1}
\\
&
+\frac{M(M-1)}{NM-1}(2N+N^2)=-(M-1)(4+N) -\frac{2}{N}\sum_{k=1}^{M-1}\frac{1}{k-\frac{1}{N}}
\\
&
+\frac{M(M-1)N(2+N)}{NM-1}.
\end{align*}

We will use the functional relation $\psi(x+1) = \psi(x) +\frac{1}{x}$ for the digamma function, which allows us to obtain, for instance,
\[
\sum_{k=1}^{M-1}\frac{1}{k-\frac{1}{N}} = \psi\left(M-\frac{1}{N}\right)-\psi\left(1-\frac{1}{N}\right).
\]
Using this we get 
\[\Sigma_2=-(M-1)(4+N)-\frac{2}{N}\left(\psi\left(M-\frac{1}{N}\right)-\psi\left(1-\frac{1}{N}\right)\right)
+\frac{M(M-1)N(2+N)}{NM-1}.
\]
We can simplify $\Sigma_1$ with the same property. Since
\[
    \sum_{k=1}^M\psi(Nk-1)=\sum_{k=1}^M\left(\psi(Nk+1)-\frac{1}{Nk-1}-\frac{1}{Nk}\right),
\]
then
\begin{align*}
    -2 & \sum_{k=0}^M \psi(Nk+1)-2\sum_{k=1}^M\psi(Nk-1)=2\gamma -4\sum_{k=1}^M \psi(Nk+1)+2\sum_{k=1}^M \frac{1}{Nk-1}+2\sum_{k=1}^M \frac{1}{Nk}
    \\
    &
    =2\gamma -4\sum_{k=1}^M \psi(Nk+1)+\frac{2}{N}\left(\psi\left(M+1-\frac{1}{N}\right)-\psi\left(1-\frac{1}{N}\right)\right)+\frac{2}{N}(\psi(M+1)+\gamma).
\end{align*}
Therefore,
\begin{align*}
& \Sigma_1+\Sigma_2=(M+2)(M+1)\psi(NM+1)+(M+1)M\psi(N(M-1)+1)
\\
&
-2(M+1)M \psi(NM+1)+(M+1)M \psi(NM-1)-M(M-1)\psi(N(M+1)-1)+2\gamma 
\\
&
+\frac{2}{N}\left(\psi\left(M+1-\frac{1}{N}\right)-\psi\left(1-\frac{1}{N}\right)\right)
+\frac{2}{N}(\psi(M+1)+\gamma)-(M-1)(4+N)
\\
&
-\frac{2}{N}\left(\psi\left(M-\frac{1}{N}\right)-\psi\left(1-\frac{1}{N}\right)\right)+\frac{M(M-1)N(2+N)}{NM-1}-4\sum_{k=1}^M \psi(Nk+1).
\end{align*}
From the relation \cite[8.365 (6)]{GR07},
\[
\sum_{k=1}^M \psi(Nk+1)= \frac{1}{N}\sum_{k=1}^M \sum_{j=1}^{N} \psi\left(k+\frac{j}{N}\right)+M\log N= \frac{1}{N} \sum_{j=1}^{N}\sum_{k=1}^M \psi\left(k+\frac{j}{N}\right)+M\log N.
\]
Summation by parts gives
\begin{align*}
    \sum_{k=1}^M & \psi\left(k+\frac{j}{N}\right)= M\psi\left(M+\frac{j}{N}\right)-\sum_{l=1}^{M-1}\left( \psi\left(l+1+\frac{j}{N}\right)-\psi\left(l+\frac{j}{N}\right)\right)l
    \\
    &
    = M\psi\left(M+\frac{j}{N}\right)-\sum_{l=1}^{M-1}\frac{l}{l+\frac{j}{N}}= M\psi\left(M+\frac{j}{N}\right)-(M-1)+\frac{j}{N}\sum_{l=1}^{M-1}\frac{1}{l+\frac{j}{N}}
    \\
    &
    = M\psi\left(M+\frac{j}{N}\right)-(M-1)+\frac{j}{N}\left(\psi\left(M+\frac{j}{N}\right)-\psi\left(1+\frac{j}{N}\right)\right)
    \\
    &
    = \left(M+\frac{j}{N}\right)\psi\left(M+\frac{j}{N}\right)-\frac{j}{N}\psi\left(1+\frac{j}{N}\right)-(M-1),
\end{align*}
for every $1\leq j \leq N$ and then
\begin{align*}
& \Sigma_1+\Sigma_2 =(M+2)(M+1)\psi(NM+1)+(M+1)M\psi(N(M-1)+1)
\\
&
-2(M+1)M \psi(NM+1)+(M+1)M \psi(NM-1)-M(M-1)\psi(N(M+1)-1)
\\
&
+2\gamma +\frac{2}{N}\left(\psi\left(M+1-\frac{1}{N}\right)-\psi\left(1-\frac{1}{N}\right)\right)+\frac{2}{N}(\psi(M+1)+\gamma)
\\
&
-(M-1)(4+N)-\frac{2}{N}\left(\psi\left(M-\frac{1}{N}\right)-\psi\left(1-\frac{1}{N}\right)\right)+\frac{M(M-1)N(2+N)}{NM-1}
\\
&
 -\frac{4}{N} \sum_{j=1}^{N}\left(M+\frac{j}{N}\right)\psi\left(M+\frac{j}{N}\right)  + \frac{4}{N} \sum_{j=1}^{N}\frac{j}{N}\psi\left(1+\frac{j}{N}\right)+4(M-1)-4M\log N
 \\
 &
 =-M(M-1)\psi(N(M+1)-1)+2(M+1)\psi(NM+1)
 \\
 &
 +(M+1)M\psi(N(M-1)+1)-\frac{M(M+1)}{NM}+\frac{(M-1) \left(\left(N^2+2 N-1\right) M-2\right)}{NM-1}
 \\
 &
 +\frac{2}{N}\psi(M+1)-N(M-1)-4M\log N-\frac{4}{N} \sum_{j=1}^{N}\left(M+\frac{j}{N}\right)\psi\left(M+\frac{j}{N}\right)
 \\
 &
  + \frac{4}{N} \sum_{j=1}^{N}\frac{j}{N}\psi\left(1+\frac{j}{N}\right) +2\gamma\left(1+\frac{1}{N}\right).
\end{align*}
Using the asymptotic expansion $\psi(z)=\log z -\frac{1}{2z}-\frac{1}{12z^2}+O\left(z^{-4}\right)$ as $z\to \infty$, we obtain
\[
\Sigma_1+\Sigma_2
=-1-\frac{3}{N}+2\log N +O\left(M^{-1}\right)+ \frac{4}{N} \sum_{j=1}^{N}\frac{j}{N}\psi\left(1+\frac{j}{N}\right)
+2\gamma\left(1+\frac{1}{N}\right).
\]
Then
\begin{align*}
& \E[E_{2}]=\frac{N^2}{2^{3}}\lim_{M\to\infty}\Big[\Sigma_1+\Sigma_2\Big]
\\
&
=\frac{N^2}{2^3}\left(-1-\frac{3}{N}+2\log N +\frac{4}{N} \sum_{j=1}^{N}\frac{j}{N}\psi\left(1+\frac{j}{N}\right)+2\gamma\left(1+\frac{1}{N}\right)\right)
\\
&
=\frac{N^2}{2^{3}}\left(-1-\frac{3}{N}+2\log N +\frac{4}{N} \underbrace{\sum_{j=1}^{N}\frac{j}{N}\frac{1}{j/N}}_{=N}+\frac{4}{N} \sum_{j=1}^{N}\frac{j}{N}\psi\left(\frac{j}{N}\right)+2\gamma\left(1+\frac{1}{N}\right)\right)
\\
&
=\frac{N}{2^{3}}\left(3N-3+2N\log N +4 \sum_{j=1}^{N}\frac{j}{N}\psi\left(\frac{j}{N}\right)+2\gamma\left(N+1\right)\right)
\\
&
=\frac{N}{2^{3}}\left(3N-3+2N\log N +4 \sum_{j=1}^{N-1}\frac{j}{N}\psi\left(\frac{j}{N}\right)+2\gamma (N-1)\right).
\end{align*}
Finally, using 
\[\sum_{j=1}^{N-1}\frac{j}{N}\psi \left(\frac{j}{N}\right)=-\frac{\gamma}{2}(N-1)-\frac{N}{2}\log N- \frac{\pi}{2}\sum_{j=1}^{N-1}\frac{j}{N}\cot\left(\frac{\pi j}{N}\right),
\]
\cite[(B.11)]{Bla15}, we get (\ref{2energy})
$$
\E[E_{2}]
=-\frac{ N\pi }{4}\sum_{j=1}^{N-1}\frac{j}{N}\cot\left(\frac{\pi j}{N}\right)+\frac{3N^2}{8}-\frac{3N}{8}.
$$

To compute $\E[E_{-2n}]$ and $\E[E_0]$, we start observing that for $r<1$ formula \eqref{expression_M} yields
\[
\E[E_{2r}]=\frac{\Gamma(1-r)N^2}{2^{2r}}\left(\frac{1}{\Gamma(2-r)}+r(1+r)\sum_{k=1}^\infty \frac{\Gamma(Nk+1)}{\Gamma(Nk+2-r)}+r(1-r)\sum_{k=1}^\infty \frac{\Gamma(Nk)}{\Gamma(Nk+2-r)}\right),
\]
since both sums are convergent in this case. Using the expression of the beta function in terms of gamma function and 
the monotone convergence theorem, we get
\begin{equation}\label{negative_r}
\E[E_{2r}]=\frac{N^2}{2^{2r}}\left(\frac{1}{1-r}+r(1+r)\int_0^1(1-t)^{-r}\frac{t^N}{1-t^N}\,dt+r\int_0^1(1-t)^{1-r}\frac{t^{N-1}}{1-t^N}\,dt\right).
\end{equation}
For $r=-n$, the energy is
\begin{equation}\label{neg2s}
\E[E_{-2n}]=2^{2n}N^2\Bigg(\frac{1}{n+1}-n(1-n)\underbrace{\int_0^1(1-t)^{n}\frac{t^N}{1-t^N}\,dt}_{I_1}-n\underbrace{\int_0^1(1-t)^{1+n}\frac{t^{N-1}}{1-t^N}\,dt}_{I_2}\Bigg).
\end{equation}
To compute $I_1$ and $I_2$ we will use the following integral representation \cite[8.361 (7)]{GR07} for the digamma function
\[
\psi(z)=\int_0^1\frac{t^{z-1}-1}{t-1}\, dt -\gamma,\qquad z>0,\]
from which we get
\begin{equation}\label{digamma_int}
\int_0^1\frac{t^a-1}{1-t^N}\,dt =\frac{1}{N}\int_0^1\frac{y^{(a+1)/N-1}-y^{1/N-1}}{1-y}\,dt
=-\frac{1}{N}\left(\psi\left(\frac{a+1}{N}\right)-\psi\left(\frac{1}{N}\right)\right),
\end{equation}
for any $a>-1$.
Then
\begin{align*}
I_1 & =\int_0^1\sum_{m=0}^n\binom{n}{m}(-1)^m\frac{t^{N+m}}{1-t^N}\,dt=\int_0^1\sum_{m=0}^n\binom{n}{m}(-1)^m\frac{t^{N+m}-1}{1-t^N}\,dt
\\
&
=-\frac{1}{N}\sum_{m=0}^n\binom{n}{m}(-1)^m\left(\psi\left(\frac{m+1}{N}+1\right)-\psi\left(\frac{1}{N}\right)\right)
\\
&
=-\frac{1}{N}\sum_{m=0}^n\binom{n}{m}(-1)^m\psi\left(\frac{m+1}{N}+1\right),
\end{align*}
where we have used $\sum_{m=0}^n\binom{n}{m}(-1)^m=0$ in the second and last equality. Applying $\psi(x+1)=\psi(x)+1/x$,
\[
I_1
=-\frac{1}{N}\sum_{m=0}^n\binom{n}{m}(-1)^m\psi\left(\frac{m+1}{N}\right)
+\sum_{m=0}^n\binom{n}{m}(-1)^{m+1}\frac{1}{m+1}
\]
and it is trivial to check that the second sum equals $-1/(n+1)$.

The integral $I_2$ can be computed in a similar way
\begin{align*}
I_2 & =\int_0^1\sum_{m=0}^{n+1}\binom{n+1}{m}(-1)^m\frac{t^{N-1+m}-1}{1-t^N}\,dt=-\frac{1}{N}\sum_{m=0}^{n+1}\binom{n+1}{m}(-1)^m\psi\left(\frac{m}{N}+1\right)
\\
&
=\frac{1}{N}\left(\gamma-\sum_{m=1}^{n+1}\binom{n+1}{m}(-1)^m\psi\left(\frac{m}{N}\right)\right)+\sum_{m=1}^{n+1}\binom{n+1}{m}(-1)^{m+1}\frac{1}{m},
\end{align*}
where the second sum is $\sum_{m=1}^{n+1}\frac{1}{m}$, as stated in \cite[0.155 (4)]{GR07}.

Finally from \eqref{neg2s} we get (\ref{evenenergies})
\begin{align*}
& \E[E_{-2n}]=2^{2n}N^2\Bigg[\frac{1}{n+1}-n(1-n)\left(-\frac{1}{N}\sum_{m=0}^n\binom{n}{m}(-1)^m\psi\left(\frac{m+1}{N}\right)-\frac{1}{n+1}\right)
\\
&
-n\left(\frac{1}{N}\left(\gamma-\sum_{m=1}^{n+1}\binom{n+1}{m}(-1)^m\psi\left(\frac{m}{N}\right)\right)+\sum_{m=1}^{n+1}\frac{1}{m}\right)\Bigg]
\\
&
=2^{2 n} N^2 \left(\frac{1}{n + 1} - \frac{n (n - 1)}{n + 1} -n \sum_{m=1}^{n+1}\frac{1}{m}\right)
\\
&
		+  2^{2 n}nN\left(-\gamma + 
    \sum_{m=1}^{n+1}\binom{n + 1}{m}(-1)^m \psi\left(\frac{m}{N}\right)\left(\frac{n - 1 }{n + 1}m + 1\right)\right).
\end{align*}

In order to compute $\E[E_0],$ i.e. formula (\ref{ourABS}) from \cite{ABS11}, we take the derivative of $\E[E_s]$ at $s=0$. Consider the continuous function
\[g(r)
     = \left\{\begin{array}{lr}
       \E[E_{2r}], & \text{for } r\neq 0,\\
        N^2-N, & \text{for } r=0,
        \end{array}\right.
  \]
where $r=0$ matches the Riesz $0$-energy, which trivially is $N^2-N$ for any configuration of points. Then
$$\E[E_0]=\frac{1}{2}g'(0).$$
Since $g'(0)$ exists, we can derive it by restricting to $r<0$
\[g'(0)=\lim_{r\to 0^{-}}\frac{g(r)-g(0)}{r},\]
where 
\[
g(r)=2^{-2r}N^2\left(\frac{1}{1-r}+r(1+r)\int_0^1(1-t)^{-r}\frac{t^N}{1-t^N}\,dt+r\int_0^1(1-t)^{1-r}\frac{t^{N-1}}{1-t^N}\,dt\right),
\] according to \eqref{negative_r}.

Then 
\begin{align*}
\lim_{r\to 0^{-}} & g(r)=N^2+N^2\lim_{r\to 0^{-}}r(1+r)\int_0^1(1-t)^{-r}\frac{t^N}{1-t^N}\,dt +N^2\lim_{r\to 0^{-}}r\int_0^1(1-t)^{1-r}\frac{t^{N-1}}{1-t^N}\,dt 
\\
&
=N^2+N^2\lim_{r\to 0^{-}}r(1+r)\int_0^1(1-t)^{-r}\frac{t^N}{1-t^N}\,dt,
\end{align*} 
because $(1-t)^{1-r}\uparrow (1-t)$ when $r\to 0^-$ and $\int_0^1(1-t)\frac{t^{N-1}}{1-t^N}\,dt<\infty$. By continuity, we also have $\lim_{r\to 0^{-}}g(r) =g(0)=N^2-N$, so we deduce that
\begin{equation}\label{limit_neg_r}
\lim_{r\to 0^{-}}r\int_0^1(1-t)^{-r}\frac{t^N}{1-t^N}\,dt=-\frac{1}{N}.
\end{equation}
Therefore,
\begin{align}\label{derg0}
\frac{g'(0)}{N^2} & =\left(1-\log 4\right)+\lim_{r\to 0^{-}}\frac{ \frac{r(1+r)}{2^{2r}}\int_0^1(1-t)^{-r}\frac{t^N}{1-t^N}\,dt+\frac{1}{N}}{r}+\lim_{r\to 0^{-}}\frac{\frac{r}{2^{2r}}\int_0^1(1-t)^{1-r}\frac{t^{N-1}}{1-t^N}\,dt}{r}
\\
&
=\left(1-\log 4\right)+\lim_{r\to 0^{-}}\frac{\frac{r}{2^{2r}}\int_0^1(1-t)^{-r}\frac{t^N}{1-t^N}\,dt+\frac{1}{N}}{r}+\underbrace{\lim_{r\to 0^{-}}r\int_0^1(1-t)^{-r}\frac{t^N}{1-t^N}\,dt}_{=-1/N\text{ by \eqref{limit_neg_r}}}\nonumber
\\
&
+\int_0^1(1-t)\frac{t^{N-1}}{1-t^N}\,dt
=\left(1-\log 4\right)
+\underbrace{\lim_{r\to 0^{-}}\frac{\frac{r}{2^{2r}}\int_0^1(1-t)^{-r}\frac{t^N}{1-t^N}\,dt+\frac{1}{N}}{r}}_{I_3}\nonumber
\\
&
-\frac{1}{N}
-\frac{1}{N}\left(\psi\left(1\right)-\psi\left(1+\frac{1}{N}\right)\right), \nonumber
\end{align}
where we have applied \eqref{digamma_int}.

It remains to compute the limit $I_3$
\begin{align*}
I_3 & =\lim_{r\to 0^{-}}\frac{\frac{r}{2^{2r}}\int_0^1(1-t)^{-r}\left(\frac{t^N}{1-t^N}-\frac{1}{N(1-t)}\right)\,dt}{r}+\frac{1}{N}\lim_{r\to 0^{-}}\frac{\frac{r}{2^{2r}}\int_0^1(1-t)^{-r}\frac{1}{1-t}\,dt+1}{r}
\\
&
=\lim_{r\to 0^{-}}\int_0^1(1-t)^{-r}\left(\frac{t^N}{1-t^N}-\frac{1}{N(1-t)}\right)\,dt
+\frac{1}{N}\lim_{r\to 0^{-}}\frac{-2^{-2r}+1}{r}
\\
&
=\underbrace{\int_0^1 \left(\frac{t^N}{1-t^N}-\frac{1}{N(1-t)}\right)\,dt}_{I_4}
+\frac{2}{N}\log 2,
\end{align*}
where the limit of the last integral is justified by monotone convergence theorem. Using \eqref{digamma_int}, we obtain
\begin{align*}
I_4 & =\frac{1}{N}\int_0^1 \frac{Nt^N-\sum_{j=0}^{N-1}t^j}{1-t^N}\,dt
=\frac{1}{N}\sum_{j=0}^{N-1}\int_0^1 \frac{t^N-t^j}{1-t^N}\,dt
\\
&
=-\frac{1}{N^2}\sum_{j=0}^{N-1}\left(\psi\left(1+\frac{1}{N}\right)-\psi\left(\frac{j+1}{N}\right)\right)
=-\frac{1}{N}\psi\left(1+\frac{1}{N}\right)+\frac{1}{N^2}\sum_{j=0}^{N-1}\psi\left(\frac{j+1}{N}\right)
\\
&
=
-\frac{1}{N}\psi\left(1+\frac{1}{N}\right)-\frac{1}{N}\left(\log N+\gamma \right),
\end{align*}
where we have used that $\sum_{j = 0}^{N - 1}\psi\left(\frac{j+1}{N}\right)=- N\log N-\gamma N.$

From \eqref{derg0} we finally get
\begin{align*}
2\E[E_0]  & =g'(0) =\left(1-\log 4\right)N^2+N^2\left(-\frac{1}{N}\psi\left(1+\frac{1}{N}\right)-\frac{1}{N}\left(\log N+\gamma \right)+\frac{2}{N}\log 2\right)
\\
&
-N
-N\left(\psi\left(1\right)-\psi\left(1+\frac{1}{N}\right)\right)
=\left(1-\log 4\right)N^2
-N\log N
-\left(1-\log 4\right)N.
\end{align*}

\end{proof}


\begin{prop}\label{aux_prop} Let $r<2$ and $m\geq 1$. Then
\begin{align}\label{prop_limit}
 & \lim_{M\to\infty}  \Biggr[2r(1+r)\sum_{k=1}^M \frac{\Gamma(Nk+1)}{\Gamma(Nk+2-r)}+2r(1-r)\sum_{k=1}^M \frac{\Gamma(Nk)}{\Gamma(Nk+2-r)} 
\\
&
-2(1+r)N^{r-1}M^r-r(N+r+Nr-r^2)N^{r-2}M^{r-1} \Biggr] \nonumber
\\
&
=
2r(1+r)\sum_{j=0}^{m-1}\frac{B_{2j}^{(r)}(\frac{r}{2})(1-r)_{2j}}{(2j)!}N^{r-1-2j}\zeta\left(1-r+2j,1+\frac{2-r}{2N}\right) \nonumber
\\
&
+2r(1-r)\sum_{j=0}^{m-1}\frac{B_{2j}^{(r-1)}(\frac{r-1}{2})(2-r)_{2j}}{(2j)!}N^{r-2-2j}\zeta\left(2-r+2j,1+\frac{1-r}{2N}\right)+ O\left(N^{r-1-2m}\right), \nonumber
\end{align}
when $N\to +\infty.$
\end{prop}


\proof
We will use the following Fields' approximation for the quotient of gamma functions, see \cite[Eq. 5.11.14]{DLMF} or \cite{Fie66}
\[\frac{\Gamma(z+a)}{\Gamma(z+b)}=\sum_{j=0}^{m-1}\frac{B_{2j}^{(2\rho)}(\rho)(b-a)_{2j}w^{a-b-2j}}{(2j)!}
+O\left(w^{a-b-2m}\right),
\]
as $w\to \infty$ with $|\arg (w+\rho)|<\pi$ where $a$ and $b$ are fixed complex numbers, $w=z+\rho$ and $2\rho = 1+a-b$. Then,
\begin{align}\label{first_sum}
& \sum_{k=1}^M  \frac{\Gamma(Nk+1)}{\Gamma(Nk+2-r)}=\sum_{k=1}^M \left(Nk+\frac{2-r}{2}\right)^{r-1}+\sum_{j=1}^{m-1}\frac{B_{2j}^{(r)}(\frac{r}{2})(1-r)_{2j}}{(2j)!}\sum_{k=1}^{M}\left(Nk+\frac{2-r}{2}\right)^{r-1-2j}
\\
&
+\sum_{k=1}^M O\left(\left(Nk+\frac{2-r}{2}\right)^{r-1-2m}\right)
=N^{r-1}\underbrace{\sum_{k=1}^M \left(k+\frac{2-r}{2N}\right)^{r-1}}_D \nonumber
\\
&
+\sum_{j=1}^{m-1}\frac{B_{2j}^{(r)}(\frac{r}{2})(1-r)_{2j}}{(2j)!}N^{r-1-2j}\underbrace{\sum_{k=0}^{M-1}\frac{1}{\left(k+1+\frac{2-r}{2N}\right)^{1-r+2j}}}_{E_j}
+ O\left(N^{r-1-2m}\right)\underbrace{\sum_{k=1}^M \frac{1}{k^{1-r+2m}}}_F, \nonumber
\end{align}
and
\begin{align}\label{second_sum}
& \sum_{k=1}^M \frac{\Gamma(Nk)}{\Gamma(Nk+2-r)}=\sum_{k=1}^M \left(Nk+\frac{1-r}{2}\right)^{r-2}
\\
&
+\sum_{j=1}^{m-1}\frac{B_{2j}^{(r-1)}(\frac{r-1}{2})(2-r)_{2j}}{(2j)!}\sum_{k=1}^{M}\left(Nk+\frac{1-r}{2}\right)^{r-2-2j}\nonumber
\\
&
+\sum_{k=1}^M O\left(\left(Nk+\frac{1-r}{2}\right)^{r-2-2m}\right)=N^{r-2}\underbrace{\sum_{k=0}^{M-1}\frac{1}{\left(k+1+\frac{1-r}{2N}\right)^{2-r}}}_G \nonumber
\\
&
+\sum_{j=1}^{m-1}\frac{B_{2j}^{(r-1)}(\frac{r-1}{2})(2-r)_{2j}}{(2j)!}N^{r-2-2j}\underbrace{\sum_{k=0}^{M-1}\frac{1}{\left(k+1+\frac{1-r}{2N}\right)^{2-r+2j}}}_{H_j}
+ O\left(N^{r-2-2m}\right)\underbrace{\sum_{k=1}^M \frac{1}{k^{2-r+2m}}}_I. \nonumber
\end{align}
To compute the limit as $M\to \infty$, observe that $E_j\to \zeta\left(1-r+2j,1+\frac{2-r}{2N}\right)$ and $H_j\to \zeta\left(2-r+2j,1+\frac{1-r}{2N}\right)$ for $j\geq 1$, since $1-r+2j, 2-r+2j > 1$. The 
sums $F$ and $I$ are convergent and $G$ can be written as (see \cite[Eq. 25.11.5]{DLMF})
\[
\sum_{k=1}^{M} \frac{1}{\left(k+\frac{1-r}{2N}\right)^{2-r}}
=\zeta\left(2-r,1+\frac{1-r}{2N}\right)-\frac{\left(M+\frac{1-r}{2N}\right)^{r-1}}{1-r}-(2-r)\int_{M-1}^\infty\frac{x-\lfloor x\rfloor}{\left(x+1+\frac{1-r}{2N}\right)^{3-r}}\,dx.
\]
The same formula holds to approximate $D$ for $r<1$
\begin{equation}\label{approxr<1}
\sum_{k=1}^{M} \frac{1}{\left(k+\frac{2-r}{2N}\right)^{1-r}}
=\zeta\left(1-r,1+\frac{2-r}{2N}\right)+\frac{\left(M+\frac{2-r}{2N}\right)^{r}}{r}-(1-r)\int_{M-1}^\infty\frac{x-\lfloor x\rfloor}{\left(x+1+\frac{2-r}{2N}\right)^{2-r}}\,dx,
\end{equation}
while if $r>1$, by the Euler-Maclaurin formula,
\begin{align}\label{approxr>1}
& \sum_{k=1}^M \left(k+\frac{2-r}{2N}\right)^{r-1}=\int_1^M \left(x+\frac{2-r}{2N}\right)^{r-1}\,dx+\frac{1}{2}\left[\left(M+\frac{2-r}{2N}\right)^{r-1}+\left(1+\frac{2-r}{2N}\right)^{r-1}\right]
\\
&
+(r-1)\int_{1}^M\frac{x-\lfloor x\rfloor - 1/2}{\left(x+\frac{2-r}{2N}\right)^{2-r}}\,dx 
=\frac{\left(M+\frac{2-r}{2N}\right)^{r}}{r}+\frac{\left(M+\frac{2-r}{2N}\right)^{r-1}}{2}-\frac{\left(1+\frac{2-r}{2N}\right)^{r}}{r} \nonumber
\\
&
+\frac{\left(1+\frac{2-r}{2N}\right)^{r-1}}{2}+(r-1)\int_{-\left(\frac{2-r}{2N}\right)}^M\frac{x-\lfloor x\rfloor - 1/2}{\left(x+\frac{2-r}{2N}\right)^{2-r}}\,dx-(r-1)\int_{-\left(\frac{2-r}{2N}\right)}^1\frac{x-\lfloor x\rfloor - 1/2}{\left(x+\frac{2-r}{2N}\right)^{2-r}}\,dx \nonumber
\\
&
=\frac{\left(M+\frac{2-r}{2N}\right)^{r}}{r}+\frac{\left(M+\frac{2-r}{2N}\right)^{r-1}}{2}-\frac{\left(1+\frac{2-r}{2N}\right)^{r}}{r}+\frac{\left(1+\frac{2-r}{2N}\right)^{r-1}}{2} \nonumber
\\
&
+(r-1)\int_{-\left(\frac{2-r}{2N}\right)}^M\frac{x-\lfloor x\rfloor - 1/2}{\left(x+\frac{2-r}{2N}\right)^{2-r}}\,dx-\left(\frac{2-r}{2N}\right)^{r-1}+\frac{\left(1+\frac{2-r}{2N}\right)^{r}}{r}-\frac{\left(1+\frac{2-r}{2N}\right)^{r-1}}{2} \nonumber
\\
&
=\frac{\left(M+\frac{2-r}{2N}\right)^{r}}{r}+\frac{\left(M+\frac{2-r}{2N}\right)^{r-1}}{2}
+(r-1)\int_{-\left(\frac{2-r}{2N}\right)}^M\frac{x-\lfloor x\rfloor - 1/2}{\left(x+\frac{2-r}{2N}\right)^{2-r}}dx-\left(\frac{2-r}{2N}\right)^{r-1}, \nonumber
\end{align}
 where the last integral converges for $1<r<2$ when $M\to +\infty$ to $\zeta\left(1-r,\frac{2-r}{2N}\right),$ see \cite[Eq. 25.11.26]{DLMF}.

Using all the previous computations, we get
\begin{align}\label{final_limit}
& \lim_{M\to\infty} \Biggr[ 2r(1+r)\sum_{k=1}^M \frac{\Gamma(Nk+1)}{\Gamma(Nk+2-r)}+2r(1-r)\sum_{k=1}^M \frac{\Gamma(Nk)}{\Gamma(Nk+2-r)}-2(1+r)N^{r-1}M^r 
\\
&
\left.
-r(N+r+Nr-r^2)N^{r-2}M^{r-1}\Biggr]=\lim_{M\to\infty}\left[2r(1+r)\left(N^{r-1}\sum_{k=1}^M \left(k+\frac{2-r}{2N}\right)^{r-1} \right.\right.\right. \nonumber
\\
&
\left.
+\sum_{j=1}^{m-1}\frac{B_{2j}^{(r)}(\frac{r}{2})(1-r)_{2j}}{(2j)!}\sum_{k=0}^{M-1}\frac{N^{r-1-2j}}{\left(k+1+\frac{2-r}{2N}\right)^{1-r+2j}}+O\left(N^{r-1-2m}\right) \sum_{k=1}^M \frac{1}{k^{1-r+2m}} \right) \nonumber
\\
&
+2r(1-r)\left(N^{r-2}\left(\zeta\left(2-r,1+\frac{1-r}{2N}\right)-\frac{\left(M+\frac{1-r}{2N}\right)^{r-1}}{1-r}\right.\right. \nonumber
\\
&
\left.
-\int_{M-1}^\infty\frac{(2-r)(x-\lfloor x\rfloor)}{\left(x+1+\frac{1-r}{2N}\right)^{3-r}}\,dx\right)+\sum_{j=1}^{m-1}\frac{B_{2j}^{(r-1)}(\frac{r-1}{2})(2-r)_{2j}}{(2j)!}
\sum_{k=0}^{M-1}\frac{N^{r-2-2j}}{\left(k+1+\frac{1-r}{2N}\right)^{2-r+2j}} \nonumber
\\
&
\left.\left.
 + O\left(N^{r-2-2m}\right)\sum_{k=1}^M \frac{1}{k^{2-r+2m}}\right)-2(1+r)N^{r-1}M^r-r(N+r+Nr-r^2)N^{r-2}M^{r-1} \right] \nonumber
\\
&
=\lim_{M\to\infty}\left[
2r(1+r)N^{r-1}\sum_{k=1}^M \left(k+\frac{2-r}{2N}\right)^{r-1}-2(1+r)N^{r-1}M^r \right. \nonumber
\\
&
\left.
-2rN^{r-2}\left(M+\frac{1-r}{2N}\right)^{r-1}-r(N+r+Nr-r^2)N^{r-2}M^{r-1}\right] \nonumber
\\
&
+2r(1-r)\sum_{j=0}^{m-1}\frac{B_{2j}^{(r-1)}(\frac{r-1}{2})(2-r)_{2j}}{(2j)!}N^{r-2-2j}\zeta\left(2-r+2j,1+\frac{1-r}{2N}\right) \nonumber
\\
&
+2r(1+r)\sum_{j=1}^{m-1}\frac{B_{2j}^{(r)}(\frac{r}{2})(1-r)_{2j}}{(2j)!}N^{r-1-2j}\zeta\left(1-r+2j,1+\frac{2-r}{2N}\right)+ O\left(N^{r-1-2m}\right) \nonumber.
\end{align}
Everything reduces to compute the limit appearing in \eqref{final_limit}. If $r<1$, using \eqref{approxr<1},
\begin{align*}
& \lim_{M\to\infty}\Biggr[2r(1+r)N^{r-1}\sum_{k=1}^M \left(k+\frac{2-r}{2N}\right)^{r-1}-2(1+r)N^{r-1}M^r-2rN^{r-2}\left(M+\frac{1-r}{2N}\right)^{r-1} 
\\
&
\left.
-r(N+r+Nr-r^2)N^{r-2}M^{r-1}\Biggr]=\lim_{M\to\infty}\Biggr[2r(1+r)N^{r-1}\left(\zeta\left(1-r,1+\frac{2-r}{2N}\right)\right.\right.
\\
&
\left.\left.
+\frac{\left(M+\frac{2-r}{2N}\right)^{r}}{r}-(1-r)\int_{M-1}^\infty\frac{x-\lfloor x\rfloor}{\left(x+1+\frac{2-r}{2N}\right)^{2-r}}\, dx\right)-2(1+r)N^{r-1}M^r \right]
\\
&
=2r(1+r)N^{r-1}\zeta\left(1-r,1+\frac{2-r}{2N}\right)
+2(1+r)N^{r-1}\lim_{M\to\infty}\left(\left(M+\frac{2-r}{2N}\right)^{r}-M^r\right)
\\
&
  =  2r(1+r)N^{r-1}\zeta\left(1-r,1+\frac{2-r}{2N}\right).
\end{align*}
If $r>1$, using \eqref{approxr>1},
\begin{align*}
    & \lim_{M\to\infty}\Bigg[2r(1+r)N^{r-1}\sum_{k=1}^M \left(k+\frac{2-r}{2N}\right)^{r-1}-2(1+r)N^{r-1}M^r-2rN^{r-2}\left(M+\frac{1-r}{2N}\right)^{r-1}
    \\
    &
-r(N+r+Nr-r^2)N^{r-2}M^{r-1}\Bigg]=\lim_{M\to\infty}\Bigg\{2r(1+r)N^{r-1}\Bigg(\frac{\left(M+\frac{2-r}{2N}\right)^{r}}{r}+\frac{\left(M+\frac{2-r}{2N}\right)^{r-1}}{2}
\\
&
+\int_{-\left(\frac{2-r}{2N}\right)}^M\frac{(r-1)(x-\lfloor x\rfloor - \frac{1}{2})}{\left(x+\frac{2-r}{2N}\right)^{2-r}}\, dx -\left(\frac{2-r}{2N}\right)^{r-1}\Bigg)-2(1+r)N^{r-1}M^r
\\
&
-2rN^{r-2}\left(M+\frac{1-r}{2N}\right)^{r-1}-r(N+r+Nr-r^2)N^{r-2}M^{r-1}\Bigg\}
\\
&
=\lim_{M\to\infty}\Bigg[2(1+r)N^{r-1}M^r\left(\left(1+\frac{2-r}{2NM}\right)^{r}-1 \right)
\\
&
+rN^{r-2}M^{r-1}\left((1+r)\left(1+\frac{2-r}{2NM}\right)^{r-1}N-2\left(1+\frac{1-r}{2NM}\right)^{r-1}-\left(N+r+Nr-r^2\right)\right)\Bigg]
\\
&
+2r(1+r)N^{r-1}\left(\zeta\left(1-r,\frac{2-r}{2N}\right)-\left(\frac{2-r}{2N}\right)^{r-1}\right)
\\
&
=\lim_{M\to\infty}\Bigg[rN^{r-2}M^{r-1}[\underbrace{(1+r)(2-r)+(1+r)N-2-(N+r+Nr-r^2)}_{=0}]+O\left(M^{r-2}\right)\Bigg]
\\
&
+2r(1+r)N^{r-1}\zeta\left(1-r,1+\frac{2-r}{2N}\right)=2r(1+r)N^{r-1}\zeta\left(1-r,1+\frac{2-r}{2N}\right),
\end{align*}
where we have used that $\zeta(s,a)-a^{-s}=\zeta(s,1+a)$.

Applying this limit on \eqref{final_limit} we get the desired result.

\qed


\section{Bounds for the minimal energy asymptotic expansion}\label{comparison_results}

We will start this section recalling some known results, and some conjectures, about the asymptotic expansion of the extremal energy $\mathcal{E}_s(N)$ attained by a set of $N$ points 
on the sphere $\mathbb S^2.$ For a more complete picture see \cite{BHS19}. 

The current knowledge about the asymptotic expansion of the minimal energy is far from complete even in $\mathbb S^2,$ but for $s\le -2$ the situation is well known. Indeed, the minimizers of the 
Riesz energy for $s<-2$ are points placed at each of the two endpoints of some diameter (for even $N$), \cite{Bjo56}, and for $s=-2,$ formula (\ref{equation-2})
shows that 
any configuration with center of mass at the origin attains the maximum $2 N^2.$ 

For $0<|s|<2,$ it is known that there exist $c,C>0$ (depending on $s$) such that
\begin{equation}\label{eq:2}
  -c N^{1+s/2}\le \mathcal E_s(N)-\frac{2^{1-s}}{2-s}N^2 \le -C N^{1+s/2},
\end{equation}
see \cite{RSZ94,Wag90,Wag92} and \cite{Bra06,AZ15} for improvements in the value of the constants leading to the bounds
\begin{equation}\label{alibounds}
\begin{aligned}
\En_s(N)-\frac{2^{1-s}}{2-s}\leq -\frac{\Gamma(1-s/2)}{2^s} N^{1+s/2},&\quad \text{if }0<s<2,\\
\En_s(N)-\frac{2^{1-s}}{2-s}\geq -\frac{\Gamma(1-s/2)}{2^s} N^{1+s/2},&\quad \text{if }-2<s<0,
\end{aligned}
\end{equation}
which were obtained with the bound given by the expected energy of random points from the spherical ensemble \cite{AZ15}.

In the boundary case $s=2$, it was shown in \cite[Proposition 3]{BHS12} that
\[-\frac{1}{4}N^2+O(N)\leq \mathcal E_2(N)-\frac{1}{4}N^2\log N \leq \frac{1}{4}N^2\log \log N+O(N^2),\]
and the upper bound was improved in \cite{AZ15} to
\begin{equation}\label{s2_bound_spherical}
\mathcal E_2(N)-\frac{1}{4}N^2\log N \leq \frac{\gamma}{4}N^2,
\end{equation}
where $\gamma$ is the Euler–Mascheroni constant.

For the logarithmic potential, it is known that there exists a constant, $C_{\log},$ such that
$$-0.0569\ldots \le C_{\log}\leq 2\log 2 
+\frac12\log\frac23+3\log\frac{\sqrt\pi}{\Gamma(1/3)}=-0.0556\ldots,$$ 
for which
\begin{equation}\label{eq:as}
\mathcal{E}_0(N)=\left( \frac{1}{2}-\log 2 \right)\,N^2-\frac12\,N\log N+C_{\log}\,N+o(N),\;\; N\to +\infty,
\end{equation}
see \cite{BS18,Lau21} and \cite{BL21} for a recent direct computation of the lower bound. The upper bound for $C_{\log}$ has been conjectured to be an 
equality by two different approaches \cite{BHS12,BS18}.

%

For $-2 < s < 4,$ $s\neq 0$,  the asymptotic expansion of the optimal Riesz $s$-energy has been conjectured in \cite{BHS12} to be, for $s\neq 2,$
\begin{equation} \label{conjecture}
  \mathcal E_s(N)=\frac{2^{1-s}}{2-s}N^2+\frac{(\sqrt{3}/2)^{s/2}\zeta_{\Lambda_2}(s)}{(4\pi)^{s/2}}N^{1+\frac{s}{2}}+o(N^{1+\frac{s}{2}}),\;\;\;\;N\to +\infty,  
\end{equation}
where $\zeta_{\Lambda_2}(s)$ is the zeta function of the hexagonal lattice, while for $s=2$ the conjectured expansion is
\begin{equation} \label{conjecture_s2}
  \mathcal E_2(N)=\frac{1}{4}N^2\log N+CN^2+O(1),\;\;\;\;N\to +\infty,  
\end{equation}
where $C = \frac{1}{4}\left[\gamma - \log(2\sqrt{3}\pi)\right]+\frac{\sqrt{3}}{4\pi}[\gamma_1(2/3)-\gamma_1(1/3)]\approx -0.08577$. Here, $\gamma_n(a)$ is the generalized Stieltjes constant in the Laurent expansion of the Hurwitz zeta function $\zeta(s,a)$ around $s=1$.

\begin{figure}
\centering
\includegraphics[width=0.8\textwidth]{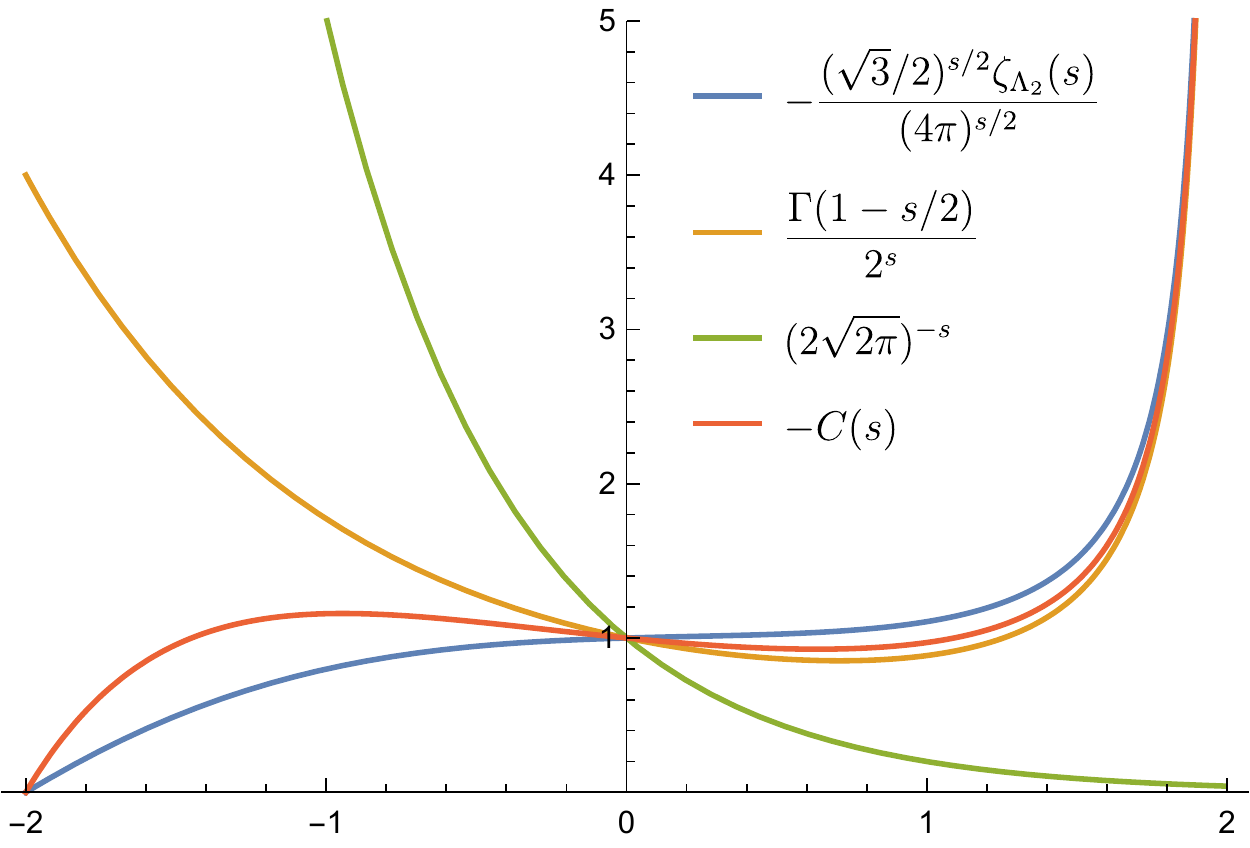}
\caption{All curves are related to $\left(\frac{2^{1-s}}{2-s}N^2-\En_s(N)\right)/N^{1+s/2}$. The 
blue curve is given by the conjectured valued for the second order constant, \eqref{conjecture}. The green and yellow curves 
corresponds to \cite{RSZ94} and \cite{AZ15}, respectively, and the red is our constant \eqref{constant1}.}
    \label{fig:grafic}
\end{figure}

It is clear that the minimal energy is always bounded by the expected energy with respect to a given random configuration. Therefore, one can bound the asymptotic expansion of the minimal energy by the 
asymptotic expansion of the expected energy. 
This idea was used in
\cite{ABS11} to get bounds for the minimal logarithmic energy using (\ref{armentano})
and
in \cite{AZ15} to get (\ref{alibounds}) and (\ref{s2_bound_spherical}). For other computations of expected energies in different settings, see \cite{BS13,BMOC16,BE18,MOC18,BE19,BF20,BDFS22,ADGMS22}. 
From our main result, Theorem \ref{asymptotics}, we obtain the asymptotic expansion \eqref{rieszasymptotics}
which is close to the conjectured expansion for
the minimal energy, see figure \ref{fig:grafic}, and we
can prove the following bounds.

\begin{cor} Let $C(s)$ be the constant in \eqref{constant1}. Then,
\begin{enumerate}[(i)]
    \item for $0<s<2$, there exists an $N_0=N_0(s)$ such that, for any $N\geq N_0$,
    \[\mathcal E_s(N)-\frac{2^{1-s}}{2-s}N^2\leq C(s) N^{1+s/2}.\]
    \item For $-2<s<0$ and a given $\epsilon>0$, there exists an $N_1=N_1(\epsilon,s)$ such that, for any $N\geq N_1$,
    \[\mathcal E_s(N)-\frac{2^{1-s}}{2-s}N^2\geq C(s)(1+\epsilon) N^{1+s/2}.\]
    \item For any $N\geq 2$,
    \begin{equation}\label{s2_bound}
    \mathcal E_2(N) - \frac{N^2\log N }{4} \leq \frac{1}{4}\left(\frac{3}{2}-\log(2\pi)+\gamma \right)N^2.
    \end{equation}
\end{enumerate}
\end{cor}


\begin{remark}
The bound (\ref{s2_bound}) improves (\ref{s2_bound_spherical}) from \cite{AZ15}, since $\frac{1}{4} \left(\frac{3}{2}-\log(2\pi)+\gamma\right)\approx 0.0598$ and $\frac{\gamma}{4} \approx 0.1443.$ In the proof we show also that
\[ \mathbb{E} [E_2]=\frac{N^2 \log N}{4} +\frac{1}{4}\left(\frac{3}{2}-\log(2\pi)+\gamma\right)N^2-\frac{N}{8}+O(1),\;\;\;\;N\to +\infty,\]
see \eqref{s2_toasymptotics}.
\end{remark}


\begin{proof}
For $0<s<2$, from 
\eqref{rieszasymptotics},
\[\frac{\mathcal E_s(N)-\frac{2^{1-s}}{2-s}N^2}{N^{1+s/2}}\leq\frac{\E[E_s]-\frac{2^{1-s}}{2-s}N^2}{N^{1+s/2}}=C(s)+\frac{s}{16}C(s-2)N^{-1}+O(N^{-2}),\quad N\to\infty.\]
Since $C(s-2)$ is negative, the last expression is bounded above by $C(s)$ for $N$ big enough.

For $-2<s<0$, using \eqref{rieszasymptotics} again,
\[\frac{\mathcal E_s(N)-\frac{2^{1-s}}{2-s}N^2}{N^{1+s/2}}\geq\frac{\E[E_s]-\frac{2^{1-s}}{2-s}N^2}{N^{1+s/2}}\xrightarrow[N\to\infty]{}C(s).\]
Therefore, given $\delta>0$, for $N$ large enough the right-hand side is bounded from below by $C(s)-\delta$. Since the constant $C(s)$ is negative, we can choose $\delta=-\epsilon C(s)$ to obtain the result.

For $s=2$, the energy is \eqref{2energy}:
\[
\E[E_{2}]
=-\frac{ N\pi }{4}\sum_{j=1}^{N-1}\frac{j}{N}\cot\left(\frac{\pi j}{N}\right)+\frac{3N^2}{8}-\frac{3N}{8}.
\]
We can rewrite the sum as
\begin{equation} \label{s2_sumexpansion}
-\sum_{j=1}^{N-1}\frac{j}{N}\cot\left(\frac{\pi j}{N}\right)
=\underbrace{\sum_{j=1}^{N-1}\left[-\frac{j}{N}\cot\left(\frac{\pi j}{N}\right)-\frac{1}{\pi(1-j/N)}\right]}_A + \underbrace{\sum_{j=1}^{N-1}\frac{1}{\pi(1-j/N)}}_B,
\end{equation}
in such a way that the term corresponding to $j=N$ in the first sum is well-defined. Let us apply the Euler-Maclaurin formula to $f(x)=g(x/N)$, with $g(x)=-x\cot(\pi x)-\frac{1}{\pi(1-x)}$:
\begin{align*}
    A
    &
    =\sum_{j=0}^{N}f(j)-f(0)-f(N)
    \\
    &
    =\int_0^N f(x) \, dx -\frac{f(0)+f(N)}{2}+\frac{B_2}{2!}[f'(N)-f'(0)]+\frac{B_4}{4!}[f^{(3)}(N)-f^{(3)}(0)]+R_N^A
    \\
    &
    =N\int_0^1 g(x) \, dx -\frac{g(0)+g(1)}{2}+\frac{1}{12N}[g'(1)-g'(0)]-\frac{1}{720N^3}[g^{(3)}(1)-g^{(3)}(0)]+R_N^A,
\end{align*}
where $B_j$ are the Bernoulli numbers and $R_N^A$ is the remainder term, that satisfies
\begin{equation}\label{remain_mclaurin}
|R_N^A|\leq \frac{2\zeta(5)}{(2\pi)^5}\int_0^N|f^{(5)}(x)|\, dx = \frac{2\zeta(5)}{(2\pi)^5N^4}\int_0^1|g^{(5)}(x)|\, dx .
\end{equation}
We get
\[
    A
    =-\frac{\log(2\pi)}{\pi}N +\frac{3}{2\pi}+\frac{\pi^2+3}{36\pi N} -\frac{\pi^4+45}{5400 \pi N^3}+R_N^A.
\]
The second sum in \eqref{s2_sumexpansion} is
\[
B=\frac{N}{\pi}\sum_{j=1}^{N-1}\frac{1}{N-j}=\frac{N}{\pi}\sum_{j=1}^{N-1}\frac{1}{j}=\frac{N}{\pi}\left(H_N-\frac{1}{N}\right),
\]
where $H_N$ is the $N$-th harmonic number. Its expansion as $N\to \infty$, see \cite{Boa77}, is
\[H_N = \log N + \gamma +\frac{1}{2N}-\frac{1}{12N^2} + R_N^H,\]
where
\begin{equation}\label{remain_harmonic}
    0<R_N^H<\frac{1}{120N^4}.
\end{equation}
With these expansions, formula \eqref{s2_sumexpansion} reads
\begin{align*}
    -\sum_{j=1}^{N-1}\frac{j}{N}\cot\left(\frac{\pi j}{N}\right)
    &
    =-\frac{\log(2\pi)}{\pi}N +\frac{3}{2\pi}+\frac{\pi^2+3}{36\pi N} -\frac{\pi^4+45}{5400 \pi N^3}+R_N^A
    \\
    &
    +\frac{N}{\pi}\left( \log N + \gamma -\frac{1}{2N}-\frac{1}{12N^2} + R_N^H\right)
    \\
    &
    =\frac{1}{\pi}\left[ N\log N +(-\log(2\pi)+\gamma)N+1+\frac{\pi^2}{36N} + NR_N^H  -\frac{\pi^4+45}{5400N^3 } +\pi R_N^A  \right].
\end{align*}
Plugging this into the formula \eqref{2energy}, we obtain
\begin{equation}\label{s2_toasymptotics}
\begin{aligned}
   \mathbb E [ E_2(x_1,\dots , x_N)]
    &
    =\frac{N}{4}\left[ N\log N +\left(\frac{3}{2}-\log(2\pi)+\gamma \right)N-\frac{1}{2}  \right.
    \\
    &
    \left. +\frac{\pi^2}{36N} + NR_N^H  -\frac{\pi^4+45}{5400 }\frac{1}{N^3} +\pi R_N^A  \right]
    \\
    &
    = \frac{N^2\log N }{4} +\frac{1}{4}\left(\frac{3}{2}-\log(2\pi)+\gamma \right)N^2-\frac{1}{8}N
    \\
    &
    +\underbrace{\frac{\pi^2}{144} + \frac{N^2R_N^H}{4}}_C  +\frac{1}{4}\Bigg(\underbrace{-\frac{\pi^4+45}{5400}\frac{1}{N^2} +\pi NR_N^A}_D  \Bigg).
\end{aligned}
\end{equation}
Finally, from \eqref{remain_harmonic}, we have
\[C \leq \frac{\pi^2}{144} + \frac{1}{480N^2} \leq \frac{\pi^2}{144} + \frac{1}{480}< 0.25 \leq \frac{N}{8}\]
for any $N\geq 2,$ and $D\leq 0$ because 
\[ \pi N|R_N^A|\leq \frac{2\pi \zeta(5)}{(2\pi)^5N^3}\int_0^1|g^{(5)}(x)|\, dx \leq  \frac{2\pi\zeta(5)}{(2\pi)^5N^3}|g^{(5)}(1)| \le  \frac{\pi^4+45}{5400}\frac{1}{N^2},\]
if $N\ge 2.$ This proves \eqref{s2_bound}.
\end{proof}

\section{Proof of theorem \ref{separation_thm}} \label{sec_separation}

\begin{proof}
Since the function $F(p,q)=\mathbf{1}_{\{|p-q|\leq t\}}$ is rotational invariant, we can apply the formula \eqref{computation2}
\begin{align*}
 & 2\E[G(t,X_N)]=\E\left[\sum_{i\neq j}\mathbf{1}_{\{|p_i-p_j|\leq t\}}\right]=\pi\int_{z  \in \CC,\,\frac{2|z|}{\sqrt{1+|z|^2}}\leq t}\rho_2(z,0)dz
\\
&
=N^2\int_0^\frac{t^2}{4-t^2} \frac{\left[\left((1+x)^N-1-Nx\right)^2(1+x)^{N-2}+\left((1+x)^N-1-Nx(1+x)^{N-1}\right)^2\right]}{[(1+x)^N-1]^3}dx,
\end{align*}
where we have applied the change of variables $r=\sqrt{x}$. As in the proof of Theorem \ref{asymptotics}, we use the identity  $\frac{1}{(x-1)^3}=\frac{1}{2}\sum_{k=2}^\infty k(k-1)x^{-(k+1)}$ for $x>1$ to get
\begin{align*}
4\E[G(t,X_N)] & =N^2\lim_{M\to\infty}\sum_{k=2}^M k(k-1)
\\
&
\Bigg[\underbrace{\int_0^\frac{t^2}{4-t^2} \left[(1+x)^{-2-Nk}+(1+x)^{-N(k+1)}\right]\left((1+x)^N-1\right)^2 dx}_{A_k} 
\\
&
-2N\underbrace{\int_0^\frac{t^2}{4-t^2} x\left[(1+x)^{-2-Nk}+(1+x)^{-1-Nk}\right]\left((1+x)^N-1\right) dx}_{B_k}
\\
&
+N^2\underbrace{\int_0^\frac{t^2}{4-t^2} x^{2}\left[(1+x)^{-2-Nk}+(1+x)^{-2-N(k-1)}\right] dx}_{C_k}\Bigg].
\end{align*}
The expression is the same than \eqref{begin_s} with $r=0$, but changing the upper limit of integration. We can take advantage of our previous computations using the following representation for the incomplete beta function
\begin{equation}\label{incbeta}
B_{s/(s+1)}(x,y)=\int_0^s \frac{t^{x-1}}{(1+t)^{x+y}}\, dt\qquad s,x,y>0,
\end{equation}
which follows from  
\[
B_{z}(x,y):= \int_0^z u^{x-1}(1-u)^{y-1}\,du\qquad 0\leq z\leq 1,\;x,y>0
\]
and the change of variables $u=t/(t+1).$

Therefore, if $s=\frac{t^2}{4-t^2}$ and $s'=\frac{s}{s+1}$, we have the analogues of \eqref{Ak}, \eqref{Bk} and \eqref{Ck}
\begin{align*}
A_k & =B_{s'}(1,N(k-2)+1)-2B_{s'}(1,N(k-1)+1)+B_{s'}(1,Nk+1)
+B_{s'}(1,N(k-1)-1)
\\
&
-2B_{s'}(1,Nk-1)+B_{s'}(1,N(k+1)-1),
\end{align*}
\[B_k=B_{s'}(2,N(k-1))-B_{s'}(2,Nk)+B_{s'}(2,N(k-1)-1)-B_{s'}(2,Nk-1),\]
and
\[C_k=B_{s'}(3,Nk-1)+B_{s'}(3,N(k-1)-1).\]

Then, changing indices as in the proof of Theorem \ref{asymptotics}, we get
\begin{align*}
4  \E[G(t,X_N)] & =N^2\lim_{M\to\infty}\sum_{k=2}^M k(k-1)\Bigg[B_{s'}(1,N(k-2)+1)-2B_{s'}(1,N(k-1)+1)
\\
&
+B_{s'}(1,Nk+1)+B_{s'}(1,N(k-1)-1)-2B_{s'}(1,Nk-1)+B_{s'}(1,N(k+1)-1)
\\
&
-2N\left(B_{s'}(2,N(k-1))-B_{s'}(2,Nk)+B_{s'}(2,N(k-1)-1)-B_{s'}(2,Nk-1)\right)
\\
&
+N^2\left(B_{s'}(3,Nk-1)+B_{s'}(3,N(k-1)-1)\right)\Bigg]=N^2\lim_{M\to\infty} \Bigg[2\sum_{k=0}^M  B_{s'}(1,Nk+1)
\\
&
+2\sum_{k=1}^M  B_{s'}(1,Nk-1)-4N\sum_{k=1}^M  kB_{s'}(2,Nk)-4N\sum_{k=1}^M  kB_{s'}(2,Nk-1)
\\
&
+2N^2\sum_{k=1}^M  k^2 B_{s'}(3,Nk-1)+g_{N,s}(M)\Bigg],
\end{align*}
where 
\begin{align*}
 g_{N,s}(M) & =(M+1)(M-2)B_{s'}(1,NM+1)-(M+1)MB_{s'}(1,N(M-1)+1)
 \\
 &
 -(M+1)MB_{s'}(1,NM-1)+M(M-1)B_{s'}(1,N(M+1)-1)
 \\
 &
 +2N(M+1)MB_{s'}(2,NM)+2N(M+1)MB_{s'}(2,NM-1)
 \\
 &
 -N^2(M+1)MB_{s'}(3,NM-1).
\end{align*}
Observe that for any $n\geq 1$,
\[B_{s'}(1,n)=\int_0^s\frac{1}{(1+t)^{n+1}}\,dt=-\frac{(1+s)^{-n}}{n}+\frac{1}{n}.
\]
Integrating by parts, one can also check that
\[B_{s'}(2,n)=-\frac{s(1+s)^{-n-1}}{n+1}-\frac{(1+s)^{-n}}{(n+1)n}+\frac{1}{(n+1)n}\]
and
\[B_{s'}(3,n)=-\frac{s^2(1+s)^{-n-2}}{n+2}-\frac{2s(1+s)^{-n-1}}{(n+2)(n+1)}-\frac{2(1+s)^{-n}}{(n+2)(n+1)n}+\frac{2}{(n+2)(n+1)n}.\]
If we replace this expressions in $g_{N,s}(M)$, we see that only the last term survives when we take the limit. For instance,
\begin{align*}
2N(M+1)MB_{s'}(2,NM) & =-\frac{2N(M+1)Ms(1+s)^{-NM-1}}{NM+1}-\frac{2N(M+1)M(1+s)^{-NM}}{(NM+1)NM}
  \\
  &
+\frac{2N(M+1)M}{(NM+1)NM},
\end{align*}
and all the terms containing the factor $(1+s)^{-NM+m}$ go to $0$ as $M\to\infty$. Therefore,
\begin{align*}
 & \lim_{M\to\infty}g_{s,N}(M)=\lim_{M\to\infty} \Bigg[\frac{(M+1)(M-2)}{NM+1}-\frac{(M+1)M}{N(M-1)+1}
-\frac{(M+1)M}{NM-1}
\\
&
+\frac{M(M-1)}{N(M+1)-1}+\frac{2N(M+1)M}{(NM+1)NM}+\frac{2N(M+1)M}{NM(NM-1)}
\\
&
-\frac{2N^2(M+1)M}{(NM+1)NM(NM-1)}\Bigg]=\lim_{M\to\infty} \left(-\frac{2}{N}+O\left(\frac{1}{M}\right)\right)=-\frac{2}{N}.
\end{align*}
On the other hand, we have
\begin{align*}
2\sum_{k=0}^M &  B_{s'}(1,Nk+1)+2\sum_{k=1}^M  B_{s'}(1,Nk-1)-4N\sum_{k=1}^M  kB_{s'}(2,Nk)-4N\sum_{k=1}^M  kB_{s'}(2,Nk-1)
\\
&
+2N^2\sum_{k=1}^M  k^2 B_{s'}(3,Nk-1)=2\sum_{k=0}^M  \frac{1-(1+s)^{-Nk-1}}{Nk+1}+2\sum_{k=1}^M \frac{1-(1+s)^{-Nk+1}}{Nk-1}
\\
&
-4N\sum_{k=1}^M k\left(-\frac{s(1+s)^{-Nk-1}}{Nk+1}+\frac{1-(1+s)^{-Nk}}{(Nk+1)Nk}\right)
\\
&
-4N\sum_{k=1}^M  k\left(-\frac{s(1+s)^{-Nk}}{Nk}+\frac{1-(1+s)^{-Nk+1}}{Nk(Nk-1)}\right)
\\
&
+2N^2\sum_{k=1}^M k^2\left(-\frac{s^2(1+s)^{-Nk-1}}{Nk+1}-\frac{2s(1+s)^{-Nk}}{(Nk+1)Nk}+\frac{2-2(1+s)^{-Nk+1}}{(Nk+1)Nk(Nk-1)}\right)
\\
&
=
2 \left(1-(1+s)^{-1}\right)+2 s\sum_{k=1}^M  (1 + s)^{- Nk-1} (2 + s -sNk).
\end{align*}
Hence,
\begin{align*}
4 & \E[G(t,X_N)]=N^2 \left(2 \left(1-(1+s)^{-1}\right)+2 s\sum_{k=1}^\infty  (1 + s)^{- Nk-1} (2 + s -sNk)-\frac{2}{N}\right)
\\
&
=N^2 \left(\frac{2s}{1+s}  -\frac{2}{N}+\frac{2s(2+s)}{(1+s)\left((1+s)^N-1\right)}-\frac{2Ns^2(1+s)^{N}}{(1+s)\left((1+s)^N-1\right)^2}\right)
\end{align*}
and with the change $s=t^2/(4-t^2)$ we get the result (\ref{the_sepa_result}).

Now we prove inequality \eqref{sep_inequality}. In terms of $s$, since $t^2=4s/(1+s)$, it reads
\[ \frac{N^2}{4} \left(-\frac{2}{N}+\frac{2s}{1+s}  +\frac{2s(2+s)}{(1+s)\left((1+s)^N-1\right)}-\frac{2Ns^2(1+s)^{N}}{(1+s)\left((1+s)^N-1\right)^2}\right)\leq \frac{N^3s^2}{8(1+s)^2},\]
or, by regrouping terms,
\[s^2 + \frac{4}{N^2}(1+s) \left(1-\frac{Ns}{(1+s)^N-1}\right)\left(1-\frac{Ns}{(1+s)^N-1}-(N-1)s\right)\geq 0.\]
Then, if we multiply by $((1+s)^N-1)^2N^2$, we have to prove
\begin{align*}
f_N(s)&:=s^2N^2((1+s)^N-1)^2
\\
&
 + 4(1+s)\left((1+s)^N-1-Ns\right) \left((1+s)^N-1-Ns-(N-1)s ((1+s)^N-1)\right) \geq 0.
\end{align*}
We expand the polynomial $(1+s)^N$ and rearrange terms in order to identify the coefficients of the polynomial $f_N$
\begin{align*}
f_N(s)&=s^2N^2\sum_{j=1}^N \binom{N}{j}s^j\sum_{k=1}^N \binom{N}{k}s^k
\\
&
+ 4(1+s) \sum_{j=2}^N \binom{N}{j}s^j \left(\sum_{k=2}^N \binom{N}{k}s^k-(N-1)s \sum_{k=1}^N \binom{N}{k}s^k\right)
\\
&
=
s^4N^2\sum_{j=0}^{N-1} \binom{N}{j+1}s^j\sum_{k=0}^{N-1} \binom{N}{k+1}s^k
\\
&
+ 4s^4 \underbrace{(1+s)\sum_{j=0}^{N-2} \binom{N}{j+2}s^j \left(\sum_{k=0}^{N-2} \binom{N}{k+2}s^k-(N-1) \sum_{k=0}^{N-1} \binom{N}{k+1}s^k\right)}_A.
\end{align*}
Using that $\binom{n}{k} = 0$ if $k>n$, the expression $A$ can be expanded in the following way
\begin{align*}
A
&
=
(1+s)\sum_{j=0}^{N-2} \binom{N}{j+2}s^j \left(\sum_{k=0}^{N-1} \binom{N}{k+2}s^k-(N-1) \sum_{k=0}^{N-1} \binom{N}{k+1}s^k\right)
\\
&
=
\left(\sum_{j=0}^{N-1} \left(\binom{N}{j+2}+\binom{N}{j+1}\right)s^j-N\right) \left(\sum_{k=0}^{N-1} \binom{N}{k+2}s^k-(N-1) \sum_{k=0}^{N-1} \binom{N}{k+1}s^k\right)
\\
&
=
\sum_{j=0}^{N-1}\binom{N}{j+2}s^j\sum_{k=0}^{N-1} \binom{N}{k+2}s^k-(N-1)\sum_{j=0}^{N-1}\binom{N}{j+2}s^j\sum_{k=0}^{N-1} \binom{N}{k+1}s^k
\\
&
+
\sum_{j=0}^{N-1}\binom{N}{j+1}s^j\sum_{k=0}^{N-1} \binom{N}{k+2}s^k
-(N-1)\sum_{j=0}^{N-1}\binom{N}{j+1}s^j\sum_{k=0}^{N-1} \binom{N}{k+1}s^k
\\
&
-N \sum_{k=0}^{N-1}\left( \binom{N}{k+2}-(N-1)  \binom{N}{k+1}\right)s^k,
\end{align*}
and then $g_N(s):=f_N(s)/s^4$ becomes
\begin{align*}
g_N(s)
&
=
(N-2)^2\sum_{j=0}^{N-1} \binom{N}{j+1}s^j\sum_{k=0}^{N-1} \binom{N}{k+1}s^k
-4(N-2)
\sum_{j=0}^{N-1}\binom{N}{j+1}s^j\sum_{k=0}^{N-1} \binom{N}{k+2}s^k
\\
&
+4\sum_{j=0}^{N-1}\binom{N}{j+2}s^j\sum_{k=0}^{N-1} \binom{N}{k+2}s^k
-4N \sum_{k=0}^{N-1}\left( \binom{N}{k+2}-(N-1)  \binom{N}{k+1}\right)s^k.
\end{align*}
Now we compute the products of sums. For instance, for the first product,
\begin{equation*}
\sum_{j=0}^{N-1}\binom{N}{j+1}s^j\sum_{k=0}^{N-1} \binom{N}{k+1}s^k=\sum_{m=0}^{2N-2}\sum_{j=0}^{m} \binom{N}{j+1} \binom{N}{m-j+1}s^m.
\end{equation*}
The same can be done with the others, yielding
\begin{align*}
g_N(s)
=
\sum_{m=0}^{2N-2}&\Bigg[(N-2)^2\underbrace{\sum_{j=0}^{m} \binom{N}{j+1} \binom{N}{m-j+1}}_B
-4(N-2)
\underbrace{\sum_{j=0}^{m} \binom{N}{j+1} \binom{N}{m-j+2} }_C
\\
&
+4\underbrace{\sum_{j=0}^{m} \binom{N}{j+2} \binom{N}{m-j+2}}_D -4N \left( \binom{N}{m+2}-(N-1)  \binom{N}{m+1} \right) \Bigg]s^m.
\end{align*}
Next, we apply Vandermonde's identity, \cite[0.156]{GR07}, to obtain the sums $B$, $C$ and $D$:
\begin{align*}
B&
=
\sum_{k=1}^{m+1} \binom{N}{k} \binom{N}{m+2-k}
=
\sum_{k=0}^{m+2} \binom{N}{k} \binom{N}{m+2-k}-2\binom{N}{m+2}
\\
&
=\binom{2N}{m+2}-2\binom{N}{m+2}
\\
C
&
=
\binom{2N}{m+3}-2\binom{N}{m+3}-N\binom{N}{m+2}
\\
D
&
=
\binom{2N}{m+4}-2\binom{N}{m+4}-2N\binom{N}{m+3}
.
\end{align*}
Then
\begin{align*}
g_N(s)
=
\sum_{m=0}^{2N-2}&\left[(N-2)^2 \left( \binom{2N}{m+2}-2\binom{N}{m+2}\right)\right.
\\
&
-4(N-2)
\left(\binom{2N}{m+3}-2\binom{N}{m+3}-N\binom{N}{m+2} \right)
\\
&
+4\left(\binom{2N}{m+4}-2\binom{N}{m+4}-2N\binom{N}{m+3}\right)\\
&
\left. 
-4N \left( \binom{N}{m+2}-(N-1)  \binom{N}{m+1} \right) \right]s^m =: \sum_{m=0}^{2N-2} c_{N,m}s^m.
\end{align*}
Remember that our goal is to check that $g_N(s)\geq 0$ for $s\geq 0$. In fact, we will see that the coefficients of this polynomial are all positive for any $N\geq 2$. To prove this, let us successively apply the identity $\binom{n}{k}=\frac{n-k+1}{k}\binom{n}{k-1}$ to get
\[
c_{N,m}
=
\frac{1}{(m+2)(m+3)(m+4)}\left(\binom{N}{m+1}h^A(N,m) +\binom{2N}{m+1}h^B(N,m) \right),
\]
where
\begin{align*}
h^A(N,m)&
=
2N^3(m^2+7m+8)+2N^2(m^3+8m^2+23m+4)-4N(m^2+m+10)+16(m+1),\\
h^B(N,m)&
=
2N^3(m^2-m-4)-N^2(m^3+3m -20)+4N(m^2-m+2)-8(m+1).
\end{align*}

Using the trivial inequality
\[2N^2(m^3+8m^2+23m)-4N(m^2+m)+16(m+1)\geq 0,\]
we have
\begin{equation}\label{bound_coefficients}
    c_{N,m}
\geq
\frac{1}{(m+2)(m+3)(m+4)}\left(\binom{N}{m+1}h^C(N,m) +\binom{2N}{m+1}h^B(N,m) \right),
\end{equation}
where
\[
h^C(N,m) = 2N^3(m^2+7m+8)+8N^2-40N.
\]

Now we check that $c_{N,m}\geq 0$, for any $0\leq m \leq 2L-2$, $L\geq 2$. Let $m\geq 3$. Then, both $h^C$ and $h^B$ are positive. Indeed, for the first one,
\[
h^C(N,m) \geq  16N^3+8N^2-40N = 8N(2N^2+N-5)\geq 0.
\]
For the second one, taking into account that we restrict to $m\leq 2N-2$, which means $N\geq (m+2)/2$, it is easy to see that $h^B(N,m)$ is increasing as a function of $N$ and therefore
\[h^B(N,m)\geq h^B(\lceil m/2 \rceil + 1,m)\geq h^B(m/2 + 1,m)= \frac{1}{4} (m-2)^2 \left(m^2+7 m+12\right)\geq 0.\]
Finally, from \eqref{bound_coefficients}, we also have
\[
(m+2)(m+3)(m+4)c_{N,m}\geq
\begin{cases}
  8 N \left(6 N^2-3 N-2\right),  & m=0, \\
  4 N \left(7 N^3-6 N^2-5 N+4\right),  & m=1, \\
  2N/3 (5 N^5- 13 N^4 + 20 N^3 - 56 N^2 + 68 N-24),  & m=2, \\
\end{cases}
\]
which are positive for any $N\geq 2$, and we are done.
\end{proof}

\bibliographystyle{unsrt}

\end{document}